\documentclass[oneside, 12pt]{amsart}

\usepackage{amsfonts}
\usepackage{amssymb}
\usepackage{latexsym}
\usepackage{graphics}
\usepackage[2emode]{psfrag}
\usepackage{mathrsfs}
\usepackage{amsthm}
\usepackage{amsmath}
\usepackage[all]{xy}
\usepackage{enumerate}
\usepackage{hyperref}

\DeclareMathOperator{\newtriangle}{\raisebox{0ex}{\scalebox{0.75}{\ensuremath{\triangledown}}}}
\DeclareMathOperator{\newblacktriangle}{\raisebox{0ex}{\scalebox{0.75}{\ensuremath{\blacktriangledown}}}}
\DeclareMathOperator{\newtriangled}{\raisebox{0.1ex}{\scalebox{0.5}{\ensuremath{\triangledown}}}}
\DeclareMathOperator{\newblacktriangled}{\raisebox{0.1ex}{\scalebox{0.5}{\ensuremath{\blacktriangledown}}}}
\DeclareMathOperator{\newdiamond}{\raisebox{0.2ex}{\scalebox{0.75}{\ensuremath{\lozenge}}}}
\DeclareMathOperator{\newblackdiamond}{\raisebox{0.2ex}{\scalebox{0.75}{\ensuremath{\blacklozenge}}}}
\DeclareMathOperator{\newdiamondd}{\raisebox{0.2ex}{\scalebox{0.5}{\ensuremath{\lozenge}}}}
\DeclareMathOperator{\newblackdiamondd}{\raisebox{0.2ex}{\scalebox{0.5}{\ensuremath{\blacklozenge}}}}

\begin{document}

\newcommand{\thlabel}[1]{\label{th:#1}}
\newcommand{\thref}[1]{Theorem~\ref{th:#1}}
\newcommand{\selabel}[1]{\label{se:#1}}
\newcommand{\seref}[1]{Section~\ref{se:#1}}
\newcommand{\lelabel}[1]{\label{le:#1}}
\newcommand{\leref}[1]{Lemma~\ref{le:#1}}
\newcommand{\prlabel}[1]{\label{pr:#1}}
\newcommand{\prref}[1]{Proposition~\ref{pr:#1}}
\newcommand{\colabel}[1]{\label{co:#1}}
\newcommand{\coref}[1]{Corollary~\ref{co:#1}}
\newcommand{\relabel}[1]{\label{re:#1}}
\newcommand{\reref}[1]{Remark~\ref{re:#1}}
\newcommand{\exlabel}[1]{\label{ex:#1}}
\newcommand{\exref}[1]{Example~\ref{ex:#1}}
\newcommand{\delabel}[1]{\label{de:#1}}
\newcommand{\deref}[1]{Definition~\ref{de:#1}}
\newcommand{\eqlabel}[1]{\label{eq:#1}}
\newcommand{\equref}[1]{(\ref{eq:#1})}

\newcommand{\norm}[1]{\| #1 \|}
\def\N{\mathbb N}
\def\Z{\mathbb Z}
\def\Q{\mathbb Q}
\def\mod{\textit{\emph{~mod~}}}
\def\R{\mathcal R}
\def\S{\mathcal S}
\def\*c{{}^*\hspace*{-1pt}{\cc}}
\def\C{\mathcal C}
\def\D{\mathcal D}
\def\J{\mathcal J}
\def\M{\mathcal M}
\def\T{\mathcal T}          

\newcommand{\Hom}{{\rm Hom}}
\newcommand{\End}{{\rm End}}
\newcommand{\Ext}{{\rm Ext}}
\newcommand{\Fun}{{\rm Fun}}
\newcommand{\Mor}{{\rm Mor}\,}
\newcommand{\Aut}{{\rm Aut}\,}
\newcommand{\Hopf}{{\rm Hopf}\,}
\newcommand{\Ann}{{\rm Ann}\,}
\newcommand{\Ker}{{\rm Ker}\,}
\newcommand{\Coker}{{\rm Coker}\,}
\newcommand{\im}{{\rm Im}\,}
\newcommand{\coim}{{\rm Coim}\,}
\newcommand{\Trace}{{\rm Trace}\,}
\newcommand{\Char}{{\rm Char}\,}
\newcommand{\Mod}{{\bf mod}}
\newcommand{\Spec}{{\rm Spec}\,}
\newcommand{\Span}{{\rm Span}\,}
\newcommand{\sgn}{{\rm sgn}\,}
\newcommand{\Id}{{\rm Id}\,}
\newcommand{\Com}{{\rm Com}\,}
\newcommand{\codim}{{\rm codim}}
\newcommand{\Mat}{{\rm Mat}}
\newcommand{\Coint}{{\rm Coint}}
\newcommand{\Incoint}{{\rm Incoint}}
\newcommand{\can}{{\sf can}}
\newcommand{\sign}{{\rm sign}}
\newcommand{\kar}{{\rm kar}}
\newcommand{\rad}{{\rm rad}}
\newcommand{\Rat}{{\rm Rat}}
\newcommand{\sd}{{\sf d}}

\def\Ab{\underline{\underline{\rm Ab}}}
\def\lan{\langle}
\def\ran{\rangle}
\def\ot{\otimes}
\def\bul{\bullet}
\def\ubul{\underline{\bullet}}
\def\tildej{\tilde{\jmath}}
\def\barj{\bar{\jmath}}

\def\id{\textrm{{\small 1}\normalsize\!\!1}}    
\def\To{{\multimap\!\to}}
\def\bigperp{{\LARGE\textrm{$\perp$}}} 
\newcommand{\QED}{\hspace{\stretch{1}}
\makebox[0mm][r]{$\Box$}\\}

\def\AA{{\mathbb A}}
\def\BB{{\mathbb B}}
\def\CC{{\mathbb C}}
\def\DD{{\mathbb D}}
\def\EE{{\mathbb E}}
\def\FF{{\mathbb F}}
\def\GG{{\mathbb G}}
\def\HH{{\mathbb H}}
\def\II{{\mathbb I}}
\def\JJ{{\mathbb J}}
\def\KK{{\mathbb K}}
\def\LL{{\mathbb L}}
\def\MM{{\mathbb M}}
\def\NN{{\mathbb N}}
\def\OO{{\mathbb O}}
\def\PP{{\mathbb P}}
\def\QQ{{\mathbb Q}}
\def\RR{{\mathbb R}}
\def\SS{{\mathbb S}}
\def\TT{{\mathbb T}}
\def\UU{{\mathbb U}}
\def\VV{{\mathbb V}}
\def\WW{{\mathbb W}}
\def\XX{{\mathbb X}}
\def\YY{{\mathbb Y}}
\def\ZZ{{\mathbb Z}}

\def\aa{{\mathfrak A}}
\def\bb{{\mathfrak B}}
\def\cc{{\mathfrak C}}
\def\dd{{\mathfrak D}}
\def\ee{{\mathfrak E}}
\def\ff{{\mathfrak F}}
\def\gg{{\mathfrak G}}
\def\hh{{\mathfrak H}}
\def\ii{{\mathfrak I}}
\def\jj{{\mathfrak J}}
\def\kk{{\mathfrak K}}
\def\ll{{\mathfrak L}}
\def\mm{{\mathfrak M}}
\def\nn{{\mathfrak N}}
\def\oo{{\mathfrak O}}
\def\pp{{\mathfrak P}}
\def\qq{{\mathfrak Q}}
\def\rr{{\mathfrak R}}
\def\ss{{\mathfrak S}}
\def\tt{{\mathfrak T}}
\def\uu{{\mathfrak U}}
\def\vv{{\mathfrak V}}
\def\ww{{\mathfrak W}}
\def\xx{{\mathfrak X}}
\def\yy{{\mathfrak Y}}
\def\zz{{\mathfrak Z}}

\def\aaa{{\mathfrak a}}
\def\bbb{{\mathfrak b}}
\def\ccc{{\mathfrak c}}
\def\ddd{{\mathfrak d}}
\def\eee{{\mathfrak e}}
\def\fff{{\mathfrak f}}
\def\ggg{{\mathfrak g}}
\def\hhh{{\mathfrak h}}
\def\iii{{\mathfrak i}}
\def\jjj{{\mathfrak j}}
\def\kkk{{\mathfrak k}}
\def\lll{{\mathfrak l}}
\def\mmm{{\mathfrak m}}
\def\nnn{{\mathfrak n}}
\def\ooo{{\mathfrak o}}
\def\ppp{{\mathfrak p}}
\def\qqq{{\mathfrak q}}
\def\rrr{{\mathfrak r}}
\def\sss{{\mathfrak s}}
\def\ttt{{\mathfrak t}}
\def\uuu{{\mathfrak u}}
\def\vvv{{\mathfrak v}}
\def\www{{\mathfrak w}}
\def\xxx{{\mathfrak x}}
\def\yyy{{\mathfrak y}}
\def\zzz{{\mathfrak z}}

\newcommand{\aA}{\mathscr{A}}
\newcommand{\bB}{\mathscr{B}}
\newcommand{\cC}{\mathscr{C}}
\newcommand{\dD}{\mathscr{D}}
\newcommand{\eE}{\mathscr{E}}
\newcommand{\fF}{\mathscr{F}}
\newcommand{\gG}{\mathscr{G}}
\newcommand{\hH}{\mathscr{H}}
\newcommand{\iI}{\mathscr{I}}
\newcommand{\jJ}{\mathscr{J}}
\newcommand{\kK}{\mathscr{K}}
\newcommand{\lL}{\mathscr{L}}
\newcommand{\mM}{\mathscr{M}}
\newcommand{\nN}{\mathscr{N}}
\newcommand{\oO}{\mathscr{O}}
\newcommand{\pP}{\mathscr{P}}
\newcommand{\qQ}{\mathscr{Q}}
\newcommand{\rR}{\mathscr{R}}
\newcommand{\sS}{\mathscr{S}}
\newcommand{\tT}{\mathscr{T}}
\newcommand{\uU}{\mathscr{U}}
\newcommand{\vV}{\mathscr{V}}
\newcommand{\wW}{\mathscr{W}}
\newcommand{\xX}{\mathscr{X}}
\newcommand{\yY}{\mathscr{Y}}
\newcommand{\zZ}{\mathscr{Z}}

\newcommand{\Aa}{\mathcal{A}}
\newcommand{\Bb}{\mathcal{B}}
\newcommand{\Cc}{\mathcal{C}}
\newcommand{\Dd}{\mathcal{D}}
\newcommand{\Ee}{\mathcal{E}}
\newcommand{\Ff}{\mathcal{F}}
\newcommand{\Gg}{\mathcal{G}}
\newcommand{\Hh}{\mathcal{H}}
\newcommand{\Ii}{\mathcal{I}}
\newcommand{\Jj}{\mathcal{J}}
\newcommand{\Kk}{\mathcal{K}}
\newcommand{\Ll}{\mathcal{L}}
\newcommand{\Mm}{\mathcal{M}}
\newcommand{\Nn}{\mathcal{N}}
\newcommand{\Oo}{\mathcal{O}}
\newcommand{\Pp}{\mathcal{P}}
\newcommand{\Qq}{\mathcal{Q}}
\newcommand{\Rr}{\mathcal{R}}
\newcommand{\Ss}{\mathcal{S}}
\newcommand{\Tt}{\mathcal{T}}
\newcommand{\Uu}{\mathcal{U}}
\newcommand{\Vv}{\mathcal{V}}
\newcommand{\Ww}{\mathcal{W}}
\newcommand{\Xx}{\mathcal{X}}
\newcommand{\Yy}{\mathcal{Y}}
\newcommand{\Zz}{\mathcal{Z}}

\def\units{{\mathbb G}_m}
\def\rightact{\hbox{$\leftharpoonup$}}
\def\leftact{\hbox{$\rightharpoonup$}}

\def\text#1{{\rm {\rm #1}}}

\def\smashco{\mathrel>\joinrel\mathrel\triangleleft}
\def\cosmash{\mathrel\triangleright\joinrel\mathrel<}

\def\ol{\overline}
\def\ul{\underline}
\def\dul#1{\underline{\underline{#1}}}
\def\Nat{\dul{\rm Nat}}
\def\Set{\dul{\rm Set}}

\renewcommand{\subjclassname}{\textup{2000} Mathematics Subject
     Classification}

\newtheorem{proposition}{Proposition}[section] 
\newtheorem{lemma}[proposition]{Lemma}
\newtheorem{corollary}[proposition]{Corollary}
\newtheorem{theorem}[proposition]{Theorem}

\theoremstyle{definition}
\newtheorem{Definition}[proposition]{Definition}
\newtheorem{example}[proposition]{Example}
\newtheorem{examples}[proposition]{Examples}

\theoremstyle{remark}
\newtheorem{remarks}[proposition]{Remarks}
\newtheorem{remark}[proposition]{Remark}

\title[Quasi-co-Frobenius Corings]{Quasi-co-Frobenius Corings as Galois comodules}
\date{\today}
\author{J. Vercruysse}
\address{Faculty of Engineering, Vrije Universiteit Brussel (VUB), B-1050 Brussels, Belgium}
\email{jvercruy@vub.ac.be}
\urladdr{homepages.vub.ac.be/\~{}jvercruy}
\thanks{The author is Postdoctoral Fellow of the Fund for Scientific Research--Flanders
(Belgium) (F.W.O.--Vlaanderen).}


\begin{abstract}
We compare several quasi-Frobenius-type properties for corings that appeared recently in literature and provide several new characterizations for each of these properties. By applying the theory of Galois comodules with a firm coinvariant ring, we can characterize a locally quasi-Frobenius (quasi-co-Frobenius) coring as a locally projective generator in its category of comodules.
\end{abstract}

\maketitle

\section*{Introduction}

The notion of what is now known as a \emph{Frobenius algebra}, and more general a \emph{quasi-Frobenius algebra}, appeared first about 100 years ago in the work of G.\ Frobenius on representation theory. 
Since the mid nineties of the previous century, there has been a revived interest in the study of Frobenius algebras, as they turned out to be an important tool in different fields, such as
Jones theory of subfactors of von Neumann algebras, topological quantum field theory, geometry of manifolds and quantum cohomology, the quantum Yang-Baxter equation 
and Yetter-Drinfeld modules. 
A modern approach to Frobenius algebras and their functorial properties can be found in \cite{CMZ}.

Quasi-Frobenius algebras (and rings) became a subject of pure algebraic research since the work of Nakayama in the late 1930's and Ikeda in the early 1950's. Since then, many equivalent characterizations of quasi-Frobenius algebras have been given. Perhaps one of the most striking properties of quasi-Frobenius algebras is the duality between their categories of left and right modules, induced by the $\Hom$-functor. In more recent years, quasi-Frobenius algebras also appear as a tool for constructing linear codes \cite{QFcodes}. A recent study of quasi-Frobenius algebras has been performed in \cite{NicYou:QF}.

Analogously, a theory of \emph{co-Frobenius} \cite{Lin:Sem} and \emph{quasi-co-Frobenius coalgebras} \cite{GTN} has been developed. In this setting, the considered categories of comodules posses exactly the dual properties of the categories modules over Frobenius and quasi-Frobenius algebras. Where (quasi-)Frobenius algebras are always finite dimensional, (quasi-)co-Frobenius coalgebras can be infinite dimensional, however, if a (quasi-)co-Frobenius coalgebra is finite dimensional, then it is the dual coalgebra of a (quasi-)Frobenius algebra. Therefore, the theory of (quasi-)co-Frobenius coalgebras can be understood as the extension of the theory of (quasi-)Frobenius algebras to the infinite dimensional case. In the theory of Hopf-algebras, the co-Frobenius property is closely related to the existence of integrals.

\emph{Corings} \cite{Sw} are defined as comonoids in the category of bimodules over a (possibly non-commutative) ring. This at first sight tiny difference with respect to the definition of usual coalgebras, which are comonoids in the category of one-sided modules over a commutative ring, has in fact far-reaching consequences. The application field of the theory of corings, that has seen a rapid development during the last decade, has shown to cover many interesting and maybe unexpected parts of mathematics. In particular, corings and comodules can be used to study ring extensions, bimodules, coalgebras and Hopf algebras, together with their relevant categories of modules (i.e.\ descent data, generalized descent data, comodules and Hopf modules, respectively). In light of this observation, it is no surprise that attempts have been made to interpret the above mentioned Frobenius properties of algebras, coalgebras and Hopf-algebras in terms of appropriate corings, this has been done in \cite{Guo}, \cite{IglNas:QF} and \cite{IV:cofrob}, in several levels of generality. In each of these papers, the quasi-Frobenius property of the $A$-coring $\cc$ is characterized by a corresponding quasi-Frobenius type property of the associated coinduction functor between the underlying category of $A$-modules and the category of $\cc$-comodules.

The aim of this paper, which can be viewed as a continuation of the work in \cite{IV:cofrob}, is twofold. 
First, after recalling the necessary definitions and proving some preliminary results in \seref{prem}, 
we investigate the relation between the different notions of `quasi-Frobenius-type' corings and functors that have been recently introduced. By applying the techniques developed in \cite{IV:cofrob}, we provide several new characterizations of these corings and functors. 
This part of the work is done in \seref{compar}. More precisely, we show that if a coring $\cc$ is both left and right QF in the sense of \cite{Guo}, it is exactly QF in the sense of \cite{IglNas:QF}, and both notions are special cases of the definition introduced in \cite{IV:cofrob} (see \thref{charQF}).

Secondly, we want to study quasi-Frobenius type corings by means of Galois theory for comodules. To this end, we develop in \seref{Galois} a theory of Galois comodules that are locally projective over the base algebra of the coring. Although this theory is a special situation of the theory developed in \cite{GTV:firm}, certain aspects are characteristic to this particular setting. In particular, where as in the general theory a firm ring is a part of the initial data of the Galois theory, in the new theory, the firm ring (in fact a ring with local units) can be reconstructed from the locally projective Galois comodule (see \thref{locunits}).

In the last section, we then combine the results of \seref{compar} and \seref{Galois} to characterize a qausi-Frobenius type coring in our main theorem as a faithfully flat Galois comodule or equivalently as (locally) projective generator in the category of left or right comodules (\thref{structureQF2}, \coref{structureQF2}). 

{\bf Notation.} Troughout this paper, we will denote the identity morphism on an object $X$ in a category $\Cc$, again by $X$. We will write $\Ab$ for the category of abelian groups. For an associative ring $A$ and a right $A$-module $M$, we denote $M^*=\Hom_A(M,A)$. If $N$ is a left $A$-module, then we write ${^*N}={_A\Hom}(M,A)$.

\section{Prelimiries}\selabel{prem}

In this section, we recall some known notions and results that will be used troughout this paper. We prove some new but rather elementary results concerning flatness over firm rings and local projectivity in comodule categories.

\subsection{Firm rings and rings with local units}

Let $R$ be a ring, not necessarily with unit. The category of all non-unital right $R$-modules is denoted by $\widetilde{\Mm}_R$. We say that $M\in\widetilde{\Mm}_R$ is \emph{firm} if and only if the multiplication map induces an isomorphism ${\sf m}:M\ot_RR \to M$, whose inverse is denoted by $\sd:M\to M\ot_RR,\ \sd(m)=m^r\ot_R r$. The category of all firm right $R$-modules and right $R$-linear maps is denoted by $\Mm_R$. The categories ${_R\widetilde{\Mm}}$ and ${_R\Mm}$ are defined in a symmetric ways.  One can easily verify that $R\in\Mm_R$ if and only if $R\in{_R\Mm}$, in which case we call $R$ a \emph{firm ring}. 
Now consider the following commutative diagram of functors, 
\[
\xymatrix{
\widetilde{\Mm}_R \ar[rr]^-{-\ot_RR} \ar@<-.5ex>[dr]_-{-\ot_RR} && \Ab 
\\
& \Mm_R \ar@<-.5ex>[ul]_J \ar@<.5ex>[ur]^U 
}
\]
Where $J$ and $U$ are forgetful functors.
Furthermore, the functor $-\ot_RR:\widetilde{\Mm}_R\to\Mm_R$ has a left adjoint $J$ and a right adjoint $\Hom_R(R,-)$ and is therefore exact. These observations lead to the following.
\begin{lemma}\lelabel{flatfirm}
Let $R$ be a firm ring. Then the following assertions are equivalent:
\begin{enumerate}[({a}i)]
\item the forgtful functor $\widetilde{U}:\widetilde{\Mm}_R\to \Ab$ is exact; 
\item the functor $J:\Mm_R\to \widetilde{\Mm}_R$ is exact;
\item the forgetful functor $U:\Mm_R\to\Ab$ is exact.
\end{enumerate}
Let $F$ be a left $R$-module 
and consider the following statements
\begin{enumerate}[(bi)]
\item $-\ot_RF:{\widetilde{\Mm}_R}\to\Ab$ is (left) exact 
\item $F$ is flat as a left $\widehat{R}$-module, where $\widehat{R}$ is the Dorroh-extension of $R$;
\item $-\ot_RF:{\Mm_R}\to \Ab$ is (left) exact.
\end{enumerate}
Then $(bi)$ is equivalent with $(bii)$ and follows by $(biii)$. If 
any of the equivalent statements of part $(a)$ hold, then the three statements are equivalent.
\end{lemma}

A left $R$-module $F$ is called \emph{flat} provided that 
statement $(bi)$ of \leref{flatfirm} holds.  
We say that $F$ is \emph{totally ({\rm or} completely) faithful} if for all $N\in\Mm_R$, the relation $N\ot_RF=0$ implies $N=0$. Finally, $F$ is termed \emph{faithfully flat} if $F$ is flat and the functor $-\ot_RF:\Mm_R\to\Ab$ reflects exact sequences. Remark the assymmetry in the definition of a faithfully flat $R$-module. By \leref{flatfirm}, this assymmetry disappears if the regular $R$-module is flat as a left $R$-module. 
Hence, the following Proposition characterizes faithfully flat $R$-modules if the regular $R$ is flat.

\begin{proposition}\prlabel{faithfullyflatfirm}
Let $R$ be a firm ring that is flat as a left $R$-module and $F$ a firm left $R$-module. The following statements are equivalent
\begin{enumerate}[(i)]
\item 
$F$ is faithfully flat as a left $R$-module; 
\item 
$F$ is flat 
and $F$ is totally faithful as a left $R$-module; 
\item 
$F$ is flat as a left $R$-module and 
for all proper right ideals $I\subset R$, we have $(R/I)\ot_RF\neq 0$ (i.e. $IF\neq F$).
\end{enumerate}
Moreover, $R$ is faithfully flat as a left $R$-module.
\end{proposition}
 
\begin{proof}
This proof is an adaption of \cite[12.17]{Wis:book}.\\
$\ul{(i)\Rightarrow(ii)}$. Take any $N\in\Mm_R$ and consider the sequence $0\to N\to 0$ in $\Mm_R$. If $N\ot_RF=0$, then the sequence $0\to N\ot_RF\to 0$ is exact in $\Ab$. By the property of part (i) we obtain that $0\to N\to 0$ also has to be exact. Consequently $N=0$.\\
$\ul{(ii)\Rightarrow(i)}$. Consider any sequence $\xymatrix{K\ar[r]^f& L\ar[r]^g & N}$ in $\Mm_R$ and suppose that the correspondig sequence 
\begin{equation}\eqlabel{eqnfg}
\xymatrix{
0\to K\ot_RF\ar[r]^-{f\ot_RF}& L\ot_RF\ar[r]^-{g\ot_RF} & N\ot_RF\to 0
}
\end{equation}
in $\Ab$ is exact. Since $R$ is flat as a left $R$-module, we know that kernels and cokernels in $\Mm_R$ can be computed in $\Ab$, hence we can consider the canonical sequence 
\[
\xymatrix{
0\to \im f\ot_RF\ar[r] & \ker g\ot_RF\ar[r]  & {\ker g/\im f}\ot_RF
}
\]
Then the exactness of the sequence \equref{eqnfg} implies that $\im f\ot_RF\cong \ker g\ot_R F$. Therefore, ${\ker g/\im f}\ot_RF=0$. Hence, we obtain from (ii) that $\ker g=\im f$, i.e. the given sequence is exact.\\
$\ul{(ii)\rightarrow(iii)}$. By part (ii), $(R/I)\ot_RF=0$ would imply $R/I=0$, i.e. $R=I$. For the last statement, consider the following exact row in $\Mm_R$,
$$\xymatrix{0\ar[r] &I\ar[r] &R \ar[r] & R/I \ar[r] &0}$$
since $F$ is flat as a left $R$-module and \leref{flatfirm}, we obtain the following commutative diagram in $\Ab$ with exact rows
\[
\xymatrix{
0\ar[r]  &I\ot_RF \ar[d]_{\mu_{I,F}} \ar[r] & R\ot_RF \ar[d]_{\mu_{R,F}} \ar[r] & (R/I)\ot_RF \ar[r] \ar[d]^\gamma &0\\
0\ar[r] &IF\ar[r] &F \ar[r] & F/IF \ar[r] &0
}
\]
Since $F$ is firm as a left $R$-module, $\mu_{R,F}$ is an isomorphism. From the properties of the diagram we immediately obtain that $\gamma$ is surjective. We can also see that $\mu_{I,F}$ is surjective. Take any $if\in IF$, then $if=\mu_{I,F}(ir\ot_R{^rf})$, where we used the firmness of $F$ and the fact that $ir\in I$ since $I$ is a right ideal in $R$. The surjectivity of $\mu_{I,F}$ implies that $\gamma$ is injective.\\
$\ul{(iii)\Rightarrow(ii)}$. Consider any $N\in\Mm_R$. We have to show that $N\ot_RF=0$ implies that $N=0$, or equivalently $N\neq 0$ implies that $N\ot_RF\neq 0$. Suppose there exists an element $0\neq n\in N$. Then we can construct a right ideal in $R$ as follows $I_n=\{r\in R~|~ nr=0\}$. Observe that $R/{I_n}$ is isomorphic to the cyclic right $R$-module $nR$. By (iii) we find 
$$0\neq R/{I_n}\ot_RF\cong nR\ot_R F \subset N\ot_R F,$$
the inclusion is a consequence of the flatness of $F$ as a left $R$-module.

The last statement follows now from \leref{flatfirm} and the fact that the regular $R$-module is always totally faithful.
\end{proof}

A ring $R$ is said to be a ring with \emph{right local units}, if for every finite set $r_1,\ldots, r_n\in R$, there exists an element $e\in R$ such that $r_i\cdot e= r_i$ for $1\le i\le n$. A ring with right local units is firm and $M\in \Mm_R$ if and only if for every finite set $m_1,\ldots, m_n\in M$ there exists an element $e\in R$ such that $m_i\cdot e=m_i$ for $1\le i\le n$. Similarly one defines rings with left and two-sided local units.

In the situation of a ring with local units, the characterization of totally faithful modules becomes easier and closer to the classical case of rings with units.

\begin{proposition}\prlabel{totallyfaithful}
Let $R$ be a ring with right local units and $F\in{_R\Mm}$. Then the following statements are equivalent
\begin{enumerate}[(i)]
\item For all $N\in\Mm_R$ and $n\in N$, $n\ot_R f=0$ for all $f\in F$ implies $n=0$;
\item For all cyclic $N\in \Mm_R$, the relation $N\ot_RF=0$ implies $N=0$;
\item For all $N\in\Mm_R$, the relation $N\ot_RF=0$ implies $N=0$.
\end{enumerate}
\end{proposition}

\begin{proof}
$\ul{(iii)\Rightarrow (ii)}$ and $\ul{(i)\Rightarrow (iii)}$ are obvious.\\
$\ul{(ii)\Rightarrow (i)}$. Consider any $N\in\Mm_R$ and $n\in N$ such that $n\ot_R f=0$ for all $f\in F$. Put $M=nR$ the cyclic right $R$-module generated by $n$, then we find that $M\ot_R F=0$, and consequently $M=0$. Let $e\in R$ be a right local unit for $n$, since $n=ne\in M$, we find $n=0$. 
\end{proof}

\subsection{Local projectivity}

Recall that an object $X$ in a Grothendieck category $\Cc$ is called finitely generated if and only if for any directed family of subobjects $\{X_i\}_{i\in I}$ of $X$ satisfying $X=\sum_{i\in I} X_i$, there exists an $i_0\in I$ such that $X=X_{i_0}$. If $X$ is a projective object in $\Cc$, this condition is known to be equivalent to the fact that $\Hom_\Cc(X,-)$ preserves coproducts. 
We call an object $X$ in $\Cc$ \emph{weakly locally projective} if any diagram with exact rows of the form
\begin{equation}\eqlabel{weaklyproj}
\xymatrix{
0 \ar[r] & E \ar[r]^i & X \ar[d]^f\\
& M \ar[r]^g & N \ar[r] & 0  
}
\end{equation}
and where $E$ is finitely generated, can be extended with a morphism $h:X\to M$ such that $g\circ h\circ i=f\circ i$. 

\begin{lemma}\lelabel{fg}
Let $F:\Cc\to\Dd$ be a covariant functor between grothendieck categories. If $F$ is exact and reflects isomorphisms, then $C$ is finitely generated in $\Cc$ if $F(C)$ is finitely generated in $\Dd$. 
\end{lemma}
\begin{proof}
Suppose that $C=\sum_{i\in I} X_i$ for a directed family of objects $\{X_i\}_{i\in I}$ in $\Cc$. Since $F$ is exact, we find that $F(C)=\sum_{i\in I} F(X_i)$ and $\{F(X_i)\}_{i\in I}$ is a directed family in $\Dd$. Then there exists an index $i_0$, such that $F(C)=F(X_{i_0})$, because $F(C)$ is finitely generated in $\Dd$. Since $F$ reflects isomorphisms, we find that $C=X_{i_0}$.
\end{proof}

If we take $\Cc=\Mm_A$, where $A$ is a ring with unit, we recover the definition of a locally projective module in the sense of \cite{ZH}, which we will call weakly locally projective right $A$-module. It is known that a right $A$-module $M$ is weakly locally projective if and only if there exist local dual bases, i.e.\ for all $m\in M$ there exist an element $e=e_i\ot_Af_i\in M\ot_A\Hom_A(M,A)$ such that $m=e\cdot m=e_if_i(m)$. If $M$ is weakly locally projective, then $M\ot_A\Hom_A(M,A)$ is a ring with left local units (see \cite{V}). More general, we call $M$ weakly $R$-locally projective if there exists a ring with local units $R$ and a ring morphism $\iota:R\to M\ot_A\Hom_A(M,A)$ such that for all $m\in M$ there exits a dual basis in the image of $\iota$. Weakly locally projective modules are always flat.

\begin{lemma}\lelabel{locprojcom}
If $R$ has right local units, then $R$ is weakly locally projective as a left $R$-module and (faithfully) flat as a right $R$-module.
\end{lemma}
\begin{proof}
For any $r\in R$, a local dual basis is given by $e\in R$ and $R\in {_R\Hom(R,R)}$, where $e$ is a right local unit for $r$. To prove flatness, take any exact sequence 
\[\xymatrix{
0\ar[rr] && M \ar[rr]^-f && N \ar[rr]^-g && P \ar[rr] && 0
}
\]
in $\widetilde{\Mm}_R$. If for any $\sum m\ot_R r\in M\ot_R R$, we have that $\sum f(m)\ot_R r=0$, then also
$\sum f(m)r=f(mr)=0$, hence $\sum mr=0$, so $\sum m\ot_Rr= mr\ot_R e=0$, where $e\in R$ is a right local unit for all $r$. Furthermore, take any $\sum n\ot_R r\in \ker(g)\ot_R R$, then $0=\sum g(n)r=\sum g(nr)$. Hence there exists an element $m\in M$ such that $f(m)=\sum nr$. Let $e\in R$ be a right local unit for all $r$. Then $\sum f(m)\ot_R e= \sum nr\ot_R e= \sum n\ot_R r$. Therefore the functor $-\ot_RR:\widetilde{\Mm}_R\to\Ab$ is left exact.
\end{proof}

\subsection{Corings and comodules}

Troughout the remaining part of this paper, $A$ will be an associative ring and $\cc$ an $A$-coring, i.e.\ a comonoid in the monoidal category ${_A\Mm_A}$ of $A$-bimodules, with coassociative comultiplication $\Delta=\Delta_\cc:\cc\to \cc\ot_A\cc,\ \Delta(c)=c_{(1)}\ot_A c_{(2)}$ and counit $\varepsilon=\varepsilon_\cc:\cc\to A$. To the coring $\cc$ we can associate its left dual ring $\*c={_A\Hom}(\cc,A)$ with unit $\varepsilon$ and where the multiplication is given by
$$f*g(c)=g(c_{(1)}f(c_{(2)})),$$
for all $f,g\in\*c$. Similarly, $\cc^*=\Hom_A(\cc,A)$ is a ring with multiplication $f*g(c)=f(g(c_{(1)})c_{(2)})$ and we have that $(\*c)^A=(\cc^*)^A={_A\Hom_A}(\cc,A)$. We denote the category of right $\cc$-comodules by $\Mm^\cc$, and for $M\in\Mm^\cc$, the coaction of $M$ by $\rho_M:M\to M\ot_A\cc,\ \rho_M(m)=m_{[0]}\ot_Am_{[1]}$.
There exists a functor $\Mm^\cc\to \Mm_{\*c}$, where for $M\in\Mm^\cc$, $m\in M$ and $f\in\*c$, we define $m\cdot f=m_{[0]}f(m_{[1]})$. If $\cc$ is weakly locally projective, then there is also a functor $\Rat:\Mm_{\*c}\to \Mm^\cc$, where for all $M\in\Mm_{\*c}$, we define $\Rat(M)$ as the largest $A$-submodule of $M$ that has a $\cc$-comodule structure such that $m\cdot f=m_{[0]}f(m_{[1]})$ for all $m\in\Rat(M)$ and $f\in \*c$.

Suppose that $\cc$ is an $A$-coring that is flat as a left $A$-module, then the category of right $\cc$-comodules is a Grothendieck category, see \cite[18.13]{BrzWis:book}, hence the definition of weakly locally projective object makes sense in $\Mm^\cc$. This is the case in particular if $\cc$ is weakly locally projective as a left $A$-module. Moreover, we have the following lemma.
\begin{lemma}
Let $\cc$ be an $A$-coring that is flat as a left $A$-module. Then for any right $\cc$-comodule the following statements hold.
\begin{enumerate}[(i)]
\item If $P$ is finitely generated as a right $A$-comodule, then $P$ is finitely generated as a right $\cc$-module;
\item If $\cc$ is weakly locally projective as a left $A$-module and $P$ is weakly locally projective as a right $\cc$-comodule, then $P$ is weakly locally projective as a right $A$-module.
\end{enumerate}
\end{lemma}
\begin{proof}
\ul{(i)}. This follows immediately from \leref{fg}.
\ul{(ii)}. Consider a diagram of the form \equref{weaklyproj} in the case $\Cc=\Mm_A$ and $X=P$. Since the functor $-\ot_A\cc:\Mm_A\to \Mm^\cc$ has a right adjoint $\Hom^\cc(\cc,-)$, it is right exact, and therefore $\xymatrix{M\ot_A\cc\ar[r]^-{f\ot_A\cc} & N\ot_A\cc \ar[r] &0}$ is an exact row in $\Mm^\cc$. Furthermore, since $\cc$ is weakly locally projective as a left $A$-module, we know by \cite[19.12]{BrzWis:book}, the finitely generated right $A$-submodule $E$ of the right $\cc$-comodule $P$ is contained in a right $\cc$-subcomodule $F$ of $P$ that is finitely genertated as a right $A$-module. Then it follows from part (i) that $F$ is also finitely generated in $\Mm^\cc$. Hence we obtain the following diagram in $\Mm^\cc$
\[
\xymatrix{
0 \ar[rr] && F \ar[rr]^-\iota && P \ar[d]^\rho\\
&&&& P\ot_A\cc \ar[d]^{f\ot_A\cc}\\
&& M\ot_A\cc \ar[rr]^-{g\ot_A\cc} && N\ot_A\cc \ar[rr] && 0  
}
\]
where $F$ contains $E$ as a right $A$-module. Since $P$ is weakly locally projective in $\Mm^\cc$, there exists a right $\cc$-comodule map $h:P\to M\ot_A\cc$ such that $(f\ot_A\cc)\circ \rho\circ\iota=(g\ot_A\cc)\circ h\circ \iota$. Applying $N\ot_A\varepsilon$ from the left on this equation, we then obtain that $f\circ i =g\circ  (P\ot_A\varepsilon)\circ h\circ i $, so $P$ is weakly locally projective as a right $A$-module.
\end{proof}

An important class of examples for corings, called \emph{comatrix corings} was introduced in \cite{EGT:comatrix}, and in a more general setup appeared in \cite{GTV:firm}. We will give an intermediate construction that will be of sufficient for our purposes. Let $A$ be a ring with unit and $R$ a ring with local units. Consider a right $A$-module $\Sigma$ that is weakly $R$-locally projective, then $\Sigma$ is in a natural way an $R$-$A$ bimodule, where the left $R$-action is induced by the map $\iota:R\to \Sigma\ot_A\Sigma^*,\ \iota(r)=x_r\ot_A \xi_r$ (we denote $\Sigma^*=\Hom_A(\Sigma,A)$), as follows
$$r\cdot y= x_r\xi_r(y),$$
for all $r\in R$ and $y\in\Sigma$. Obviously, $\Sigma\in {_R\Mm_A}$, i.e. $R$ has acts with local units on $\Sigma$, where the local units are the inverse images of the local dual bases under $\iota$. 
\begin{theorem}
With notation as above, the $A$-bimodule $\dd=\Sigma^*\ot_R\Sigma$ is an $A$-coring with comultiplication
$$\Delta(\xi\ot_R x)=\xi\ot_R e_i\ot _A f_i \ot_R x,$$
where $e_i\ot _Rf_i\in \im\iota\subset \Sigma\ot _A\Sigma^*$ is a local dual basis for $x$ and counit 
$$\varepsilon(\xi\ot_R x)=\xi(x).$$
\end{theorem}

Now suppose that $\cc$ is an $A$-coring and $\Sigma\in\Mm^\cc$. Let furthermore $R$ be a firm ring and $\Sigma$ an $R$-$\cc$ bicomodule, i.e. there is a ringmorphism $R\to \End^\cc(\Sigma)$ and $R$ is a firm left $R$-module with $R$ action induced by this ringmorphism. Then both the functor $-\ot_A\cc$ and the functor $\Hom_A(\Sigma,-)\ot_R\Sigma$ incude a comonad on $\Mm_A$. We call $\Sigma$ an \emph{$R$-$\cc$ comonadic-Galois} comodule if and only if the following natural transformation is an isomorphism for all $M\in\Mm_A$
$$\can_M:\Hom_A(\Sigma,M)\ot_R\Sigma\to M\ot_A\cc,\ \can_M(\phi\ot_R x)=\phi(x_{[0]})\ot_Ax_{[1]}.$$

If $\Sigma$ is weakly $R$-locally projective as a right $A$-module, then $\can_M$ will be an isomorphism for all $M$ if and only if $\can_A$ is an isomorphism.
In fact, $\can_A$ is a (canonical) $A$-coring morphism
\begin{equation}
\can=\can_A:\dd=\Sigma^*\ot_R\Sigma\to \cc,\ \can(\xi\ot_R x)=\xi(x_{[0]})x_{[1]}.
\end{equation}
In this situation, we call $\Sigma$ a \emph{$R$-$\cc$ Galois comodule}.

\section{Frobenius corings and their generalizations}\selabel{compar}

\subsection{The Morita context of two objects}

Let $k$ be a commutative ring and $\Cc$ a $k$-linear category, and take two objects $A$ and $B$ in $\Cc$. Recall that we can associate the following Morita context (see \cite[Remarks p 389, Examples 1.2]{Mul:Morita}) to $A$ and $B$. 
\begin{equation}\eqlabel{contextAB}
\NN(A,B)=(\End_\Cc(A),\End_\Cc(B),\Hom_\Cc(B,A),\Hom_\Cc(A,B),\circ,\bul),
\end{equation}
where both Morita maps $\circ$ and $\bul$ are defined by the composition in $\Cc$. Observe that any Morita context can be interpreted in this way.

Obviously, the objects $A$ and $B$ are isomorphic in $\Cc$ if and only if there exist an element $\barj\in \Hom_\Cc(B,A)$  and $j\in\Hom_\Cc(A,B)$ such that $j\bul \barj=B$ and $\barj\circ j=A$. Therefore, we will call a pair $(j,\barj)$ with the above property a pair of \emph{invertible elements} for a Morita context $\NN(A,B)$. If such a pair of invertible elements exists, then the Morita maps of the context are clearly surjective (i.e.\ the context is \emph{strict}).

Suppose now that $\Cc$ is a category with coproducts. Objects $A,B\in \Cc$ are called \emph{similar} if and only if there exists natural numbers $n$ and $m$, and split monomorphisms $a:A\to B^n$ and $b:B\to A^n$. Similar objects in a $k$-linear category can be characterized by means of the Morita context \equref{contextAB}, as follows from the following lemma.
\begin{lemma}\lelabel{Moritasim}
Consider the Morita context \equref{contextAB} associated to two objects $A$ and $B$ in the $k$-linear category $\Cc$ with coproducts.
\begin{enumerate}[(i)]
\item The Morita map $\circ$ is surjective if and only if $A$ is direct summand of $B$;
\item the Morita map $\bul$ is surjective if and only if $B$ is a direct summand of $A$;
\item the Morita context $\NN(A,B)$ is strict if and only if $A$ and $B$ are similar in $\Cc$.
\end{enumerate}
\end{lemma}
\begin{proof}
We only proof statement (i). The map $\circ$ is surjective if and only if there exist maps $a_i:B\to A$ and $b_i:A\to B$ such that $\sum_i a_i\circ b_i=A$. Using the universal property of the finite (co)product, this is equivalent with existence of a map $a:A\to B^n$ and $b:B^n\to A$ such that $a\circ b=A$, i.e.\ $A$ is a direct summand of $B$.
\end{proof}

We say that there exists a \emph{locally left invertible element} for $\circ$, if there exists an element $j\in\Hom_\Cc(A,B)$ such that for every finitely generated as $k$-submodule $F\subset A$, we can can find an element $\barj_F\in\Hom_\Cc(B,A)$ such that $\barj_F\circ j|_F = F$. The map $\circ$ is called \emph{locally surjective} if for every finitely generated as $k$-submodule $F\subset A$, we can find a finite number of elements $j^\ell_F\in\Hom_\Cc(A,B)$ and $j^\ell_F\in\Hom_\Cc(B,A)$ such that $ \sum_\ell\barj^\ell_F\circ j^\ell_F (a) = a$, for all $a\in F$. Clearly, if $(j,\barj)$ is a pair of invertible elements, then $j$ is locally left invertible for $\circ$ and if $j$ is locally left invertible for $\circ$ then $\circ$ is locally surjective.

\subsection{Adjoint functors and beyond}

Troughout this section, we will suppose that all categories are $k$-linear Grothendieck categories and that all functors preserve colimits.
For any two such functors $F,G:\Cc\to \Dd$, we know that $\Nat(F,G)$ is a set, hence a $k$-module (see \cite[Lemma 5.1]{CIGTN:Frob}). However, most results given here can be transferred to more general settings, considering suitable classes of functors.

Recall from \cite{IV:cofrob} that for any two functors $F:\Cc\to \Dd$ and $G:\Dd\to \Cc$
we can construct the following Morita context 
\begin{equation}\eqlabel{contextFG}
\MM(F,G)=(\Nat(G,G),\Nat(F,F)^{\rm op},\Nat(\Dd,FG),\Nat(GF,\Cc),\newdiamond,\newblackdiamond).
\end{equation}
where the connecting maps are given by the following formulas, 
\begin{eqnarray*}
(\alpha\newdiamondd\beta)_D=\beta_{GD}\circ G\alpha_{D} &{\rm and}&
(\beta\newblackdiamondd\alpha)_C=F\beta_{C}\circ \alpha_{FC},
\end{eqnarray*}
writing $\alpha\in\Nat(\Dd,FG)$, $\beta\in\Nat(GF,\Cc)$, $C\in\Cc$ and $D\in\Dd$.
It was proven in \cite[Theorem 3.4]{IV:cofrob} that $(G,F)$ is an adjoint pair if and only if there exists a pair of invertible elements for the Morita context $\MM(F,G)$, i.e.\ if and only if we can find elements $\zeta\in\Nat(\Dd,FG)$ and $\varepsilon\in\Nat(GF,\Cc)$ such that 
\begin{eqnarray}
(\zeta\newdiamondd\varepsilon)_D&=& \varepsilon_{GD} \circ G\zeta_D = GD;\eqlabel{triang1}\\
(\varepsilon\newblackdiamondd\zeta)_C&=& F\varepsilon_C \circ \zeta_{FC}= FC,\eqlabel{triang2}
\end{eqnarray}
for all $C\in\Cc$ and $D\in\Dd$. 

In \cite{Guo}, the functor $G$ was called a \emph{left quasi-adjoint} for $F$ if for some integer $n$, there exist natural transformations $\eta:\id_\Dd\to \coprod_{i=1}^n FG$ and $\zeta:\coprod_{i=1}^n GF\to \id_\Cc$ such that 
\begin{equation}
\zeta_{GD}\circ G\eta_D= GD,
\end{equation}
for all $D\in\Dd$. 
We will call $F$ a \emph{right quasi-adjoint} for $G$ if there exist similar $\eta$ and $\zeta$ that satisfy
\begin{equation}
F\zeta_C \circ \eta_FC= FD,
\end{equation}
for all $C\in\Cc$. If at the same time $F$ is a right quasi-adjoint for $G$ and $G$ is a left quasi-adjoint for $F$, then we call $(G,F)$ a \emph{quasi-adjoint pair}. 
The following result is immediate.
\begin{lemma}\lelabel{adjointGuo}
With notation as above, the following statements hold.
\begin{enumerate}[(i)]
\item The functor $G$ is a left quasi-adjoint for $F$ if and only if the Morita map $\newdiamond$ of the context \equref{contextFG} is surjective;
\item the functor $F$ is a right quasi-adjoint for $G$ if and only if the Morita map $\newblackdiamond$ of the context \equref{contextFG} is surjective;
\item $(F,G)$ is a quasi-adjoint pair if and only if the Morita context \equref{contextFG} is strict.
\end{enumerate}
\end{lemma}
If $(F,G)$ is an adjoint pair and $F$ is moreover a right quasi-adjoint for $G$, then we call $(F,G)$ a \emph{right quasi-Frobenius} pair. Similarly, an adjoint pair $(F,G)$ is called a \emph{left quasi-Frobenius pair} if $G$ is a left quasi-adjoint for $F$ and $(F,G)$ said to be a \emph{quasi-Frobenius pair} if it is at the same time a left and right quasi-Frobenius pair.

Two functors $L,R:\Aa\to \Bb$ are called \emph{similar} if and only if there exist natural transformations $\phi :L \to \coprod_{i=1}^nR$, $\psi:\coprod_{i=1}^nR\to L$, $\phi':R\to \coprod_{i=1}^nL$ and $\psi':\coprod_{i=1}^nL\to R$ such that $\psi\circ\phi=L$ and $\psi\circ\phi=R$. In \cite {IglNas:QF}, a functor $F:\Bb\to \Aa$ is then called a \emph{quasi-Frobenius functor} if and only if $F$ has a left adjoint $L$ and a right adjoint $R$ such that $L$ and $R$ are similar, in this situation $(L,F,R)$ is said to be a \emph{quasi-Frobenius} triple.

\begin{lemma}\lelabel{adjointalgebra}
Let $F:\Cc\to \Dd$ be a functor with a right adjoint $R$. Then there is an isomorphism of $k$-algebras
$$\Phi:\Nat(F,F)\cong \Nat(R,R)^{\rm op}$$
\end{lemma}
\begin{proof}
The stated isomorphism of sets can be obtained from the bijective correspondence between so-called mates of an adjuction, see \cite{KelStr:2cat}. Let us give the explicit form of the isomorphism for sake of compleetness.
Denote the unit and the counit of the adjunction respectively by $\lambda:\id_\Cc\to RF$ and $\kappa:FR\to \id_\Dd$. Then we define for all $\alpha\in\Nat(F,F)$ the natural transformation $$\Phi(\alpha)=R\kappa \circ R\alpha R\circ \lambda R\in \Nat(R,R).$$
By means of naturality and the properties of the unit and the counit, it can be checked that the inverse of $\Phi$ is given by
$$\Phi^{-1}(\beta)=\xi F\circ F\beta F\circ F\lambda,$$
for all $\beta\in\Nat(R,R)$. Finally, let us check that $\Phi$ is an algebra map. Take $\alpha$ and $\alpha'$ in $\Nat(F,F)$, and consider the following diagram,
\[
\xymatrix{
R \ar[rr]^-{\lambda R} \ar[d]_{\lambda R} && RFR \ar[rr]^-{R\alpha R} && RFR \ar[rr]^-{R\kappa} && R \ar[d]^-{\lambda R} \\
RFR \ar[rr]^-{RF\lambda R} \ar[d]_{R\alpha' R} && RFRFR \ar[rr]^-{RFR\alpha R} && RFRFR \ar[rr]^-{RFR\kappa} && RFR \ar[d]^{R\alpha' R}\\
RFR \ar[rr]^-{RF\lambda R} \ar@{=}[rrd] && RFRFR \ar[rr]^-{RFR\alpha R} \ar[d]^{R\kappa FR} && RFRFR \ar[rr]^-{RFR\kappa} && RFR\ar[d]^{R\kappa}\\
&& RFR \ar[rr]_-{R\alpha R} && RFR \ar[rr]_-{R\kappa} && R\ .
}
\]
In this diagram, the upper square commutes by the naturality of $\lambda$, the middle square commutes by the naturality of $\alpha'$ and the lower square commutes by the naturality of $\kappa$. The triangle commutes by the adjunction property of $F$ and $R$. Hence the diagram is commutative, which means exactly that $\Phi(\alpha\circ \alpha')=\Phi(\alpha')\circ\Phi(\alpha)$.
\end{proof}

Before stating our next result, recall that to any Morita context, we can associate in a canonical way two other Morita contexts without changing any of the relevant information.
The \emph{opposite} of a Morita context $\NN=(R,S,P,Q,\mu,\tau)$ is the Morita context $\NN^{\rm op}=(R^{\rm op},S^{\rm op},Q,P,\mu^{\rm op},\tau^{\rm op})$, where $\mu^{\rm op}(q\ot_{S^{op}} p)=\mu(p\ot_S q)$ and $\tau^{\rm op}(p\ot_{R^{op}} q)=\tau(q\ot_R p)$. 
The \emph{twisted} of a Morita context $\NN=(R,S,P,Q,\mu,\tau)$ is the Morita context $\NN^{\rm t}=(S,R,Q,P,\tau,\mu)$.

\begin{lemma}\lelabel{isoMoritas}
Let $F:\Cc\to\Dd$ be a functor with left adjoint $L$ and right adjoint $R$. Then we consider the Morita context $\NN(L,R)$ constructed as \equref{contextAB} and the contexts $\MM(L,F)$ and $\MM(F,R)$, constructed as \equref{contextFG}.
Then there are isomorphisms of Morita contexts $\NN^{\rm top}(L,R)\cong\MM(L,F)\cong\MM^{\rm op}(F,R)$.
\end{lemma}
\begin{proof}
The needed algebra isomorphisms follow from \leref{adjointalgebra}.

From the adjunctions $(L,F)$ and $(F,R)$, we can immediately deduce that 
$$\Nat(FL,\id_{\Dd})\cong \Nat(L,R)\cong\Nat(\id_{\Dd},FR).$$

Now denote the unit and counit of the adjunction $(F,R)$ respectively by $\lambda:\id_\Cc\to RF$ and $\kappa:FR\to \id_\Dd$. 
Take any $\alpha\in\Nat(R,L)$, and define $\alpha'\in\Nat(\id_{\Cc},LF)$ as 
$$\alpha'_C=\alpha_{F C}\circ\lambda_C.$$
Conversely, for any $\beta\in\Nat(\id_{\Cc},LF)$ we define $\beta'\in\Nat(R,L)$ by
$$\beta'_D=L\kappa_D\circ\beta_{R D}.$$
If we compute $\alpha''$, we find $\alpha''_D=L\kappa_D\circ\alpha_{FR D}\circ\lambda_{R D}$. By the naturality of $\alpha$, we know that $L\kappa_D\circ\alpha_{FR D}=\alpha_D\circ R\kappa_D$. 
Since $(F,R)$ is an adjoint pair, we obtain $R\kappa_D\circ\lambda_{R D}=R D$. Combining both identities, we find that $\alpha''_D=\alpha_D$. Similarly, we find $\beta''=\beta$, making use of the naturality of $\beta$. Hence $\Nat(\id_{\Cc},LF)\cong\Nat(R,L)$.

By a similar computation, based on the adjunction $(L,F)$, we find $\Nat(R,L)\cong\Nat(RF,\id_{\Cc})$. 

We leave it to the reader to verify that the obtained isomorphisms are bimodule maps and that they preserve the Morita maps.
\end{proof}

\begin{theorem}
Let $F:\Cc\to\Dd$ be a functor with left adjoint $L$ and right adjoint $R$. Then the following statements are equivalent
\begin{enumerate}[(i)]
\item $R$ is a left quasi-adjoint of $F$;
\item $R$ is a direct summand of $L$ (in the category of $k$-linear functors and natural transformations);
\item $F$ is a left quasi-adjoint of $L$.
\end{enumerate}
\end{theorem}
\begin{proof}
By \leref{adjointGuo}, (i) is equivalent with the fact that the map $\newdiamond$ of the context $\MM(F,R)$ is surjective. By \leref{isoMoritas}, we this is now equivalent with the fact that $\newdiamond$ is surjective in $\MM(L,F)$ and $\bul$ is surjective in $\NN(L,R)$. Again by \leref{adjointGuo}, this first statement means exactly (iii) and by \leref{Moritasim} the later statement is equivalent with (ii).
\end{proof}

A double application of the previous Theorem now gives immediately the following corresponde between the notions of quasi-Frobenius-type functors defined in \cite{Guo} and \cite{IglNas:QF}.

\begin{corollary}
Let $F:\Cc\to\Dd$ be a functor with left adjoint $L$ and right adjoint $R$. Then the following statements are equivalent
\begin{enumerate}[(i)]
\item $(F,R)$ is a quasi-Frobenius pair;
\item $(L,F,R)$ is a quasi-Frobenius triple;
\item $(L,F)$ is a quasi-Frobenius pair.
\end{enumerate}
\end{corollary}

\subsection{(Locally) quasi-Frobenius corings}

In this section we discuss the varions on Frobenius corings that were recently introduced in \cite{IglNas:QF}, \cite{Guo} and \cite{IV:cofrob}. We show how they are related to the notions of the previous section and provide several new characterizations.

Recall that $\cc$ is said to be \emph{Frobenius} if $\cc$ and $\*c$ are isomorphic as $A\hbox{-}\*c$ bimodules. This notion is left-right symmetric: $\cc$ is Frobenius if and only if $\cc$ and $\cc^*$ are isomorphic as $\cc^*$-$A$ bimodules.

An $A$-coring $\cc$ is \emph{left quasi-Frobenius} if $\cc$ is a direct summand of $\*c^n$ in the category ${_A\Mm_{\*c}}$, where $n$ is a finite integer. Similarly, we say that $\cc$ is right quasi-Frobenius if $\cc$ is a direct summand of a number of copies of $\cc^*$ in the category ${_{\cc^*}\Mm_A}$.

We say that an $A$-coring $\cc$ is \emph{quasi-Frobenius} if $\cc$ and $\*c$ are similar as $A$-$\*c$ bimodules, i.e.\ $\cc$ is a direct summand of $\*c^n$ and $\*c$ is a direct summand of $\cc^m$ for certain natural numbers $n$ and $m$. 

\begin{remark}
The definition of a quasi-Frobenius coring in \cite{IglNas:QF} was given in terms of a quasi-Frobenius functor, and shown to be equivalent with the definition stated here in \cite[Theorem 6.5]{IglNas:QF} (see also \thref{charQF} below).

The original definition in \cite{Guo} states that an $A$-coring is left quasi-Frobenius if and only if there exists maps $\pi_i\in{^\cc\Hom^\cc}(\cc\ot_A\cc,\cc)$ and elements $z_i\in \cc^A$, for $i=1,\ldots, n$, such that $\sum_i \pi_i(c\ot z_i)= c$ for all $c\in \cc$. \cite[Theorem 4.2]{Guo} claims that this definition is equivalent with the fact that $\cc$ is finitely generated and projective as a left $A$-module and $\cc$ is a direct summand of $\*c$ (as $A$-$\*c$ bimodules). However, the proof of this theorem states that there exist maps $\alpha_i:\*c\to \cc$ and $\delta_i:\cc\to \*c$ such that $\sum_i\delta_i\circ\alpha_i=\*c$. This means that $\*c$ is a direct summand of $\cc$, and not conversely, therefore the statement of \cite[Theorem 4.2]{Guo} is not correct. In fact, it means that $\cc$ is \emph{right} quasi-Frobenius in our terminology (see also \thref{charlQF} below). Since left quasi-Frobenius ring extensions in the sense of M\"uller \cite{Mul:Morita} are in correspondence with corings such that $\cc$ is a direct summand of $\*c^n$ (see \cite[Proposition 4.3]{Guo}), we believe the terminology as introduced above is the correct one.
\end{remark}

From the observations made in the previous sections, it is now clear that the quasi-Frobenius properties of an $A$-coring $\cc$ will be closely connected to the properties of the Morita context 
associated to $\cc$ and $\*c$ in ${_A\Mm_{\*c}}$,
\begin{equation}\eqlabel{contextC*C}
\NN(\cc,\*c)=({_A\End_{\*c}}(\cc),{_A\End_{\*c}}(\*c),{_A\Hom_{\*c}}(\*c,\cc),{_A\Hom_{\*c}}(\cc,\*c),\circ,\bul).
\end{equation}
The first result is immediate.
\begin{theorem}
An $A$-coring $\cc$ is Frobenius if and only if there exists a pair of invertible elements for the Morita context $\MM(\cc,\*c)$.
\end{theorem}

Recall that for any $A$-coring $\cc$, the induction functor $\Gg=-\ot_A\cc:\Mm_A\to\Mm^\cc$ has both a left adjoint, being the forgetful functor $\Ff:\Mm^\cc\to\Mm_A$, and a right adjoint, given by the Hom-functor $\Hh=\Hom^\cc(\cc,-):\Mm^\cc\to\Mm_A$. Similarly, we denote $\Gg':\cc\ot_A-:{_A\Mm}\to {^\cc\Mm}$, with left adjoint $\Ff':{^\cc\Mm}\to {_A\Mm}$ and right adjoint ${^\cc\Hom}(\cc,-):{^\cc\Mm}\to{_A\Mm}$.

\begin{theorem}\thlabel{charlQF}
Let $\cc$ be an $A$-coring and use notation as above. The following statements are equivalent.
\begin{enumerate}[(i)]
\item[(0)] The map $\circ$ of the Morita context \equref{contextC*C} is surjective; 
\item $\cc$ is left quasi-Frobenius;
\item $\cc^*$ is a direct summand of $\cc$ in the category ${_{\cc^*}\Mm_A}$;
\item $(\Ff,\Gg)$ is a left quasi-Frobenius pair;
\item $(\Gg,\Hh)$ is a left quasi-Frobenius pair;
\item $(\Ff',\Gg')$ is a right quasi-Frobenius pair;
\item $(\Gg',\Hh')$ is a right quasi-Frobenius pair;
\item the functor $\Hh$ is a direct summand of the functor $\Ff$;
\item the functor $\Ff'$ is a direct summand of the functor $\Hh'$;
\item the functor ${_A\Hom_{\*c}}(-,\cc):{_A\Mm_{\*c}}\to{_{\*c^*}\Mm}$ is a direct summand of the functor ${_A\Hom_{\*c}}(-,\*c)$;
\item the functor ${_A\Hom_{\*c}}(\cc,-):{_A\Mm_{\*c}}\to{\Mm_{\*c^*}}$ is a direct summand of the functor ${_A\Hom_{\*c}}(\*c,-)$; 
\item the functor ${_{\cc^*}\Hom_A}(-,\cc^*):{_{\cc^*}\Mm_A}\to\Mm_{\*c^*}$ is a direct summand of the functor ${_{\cc^*}\Hom_A}(-,\cc):$;
\item the functor ${_{\cc^*}\Hom_A}(\cc^*,-):{_{\cc^*}\Mm_A}\to{_{\*c^*}\Mm}$ is a direct summand of the functor ${_{\cc^*}\Hom_A}(\cc,-)$ are similar; 
\end{enumerate}
If these statements are satisfied, than $\cc$ is finitely generated and projective as a right $A$-module.
\end{theorem}

\begin{proof}
The equivalence between the first two statements follows by \leref{Moritasim}.
In \cite[Theorem 5.10]{IV:cofrob} it was proven that this Morita context is isomorphic to several contexts of natural transformations. Interpreting the surjectivity of the corresponding Morita map in each of these Morita contexts as in \leref{Moritasim} and \leref{adjointGuo} leads directly to the equivalent statements.

The last statement follows by \cite[Corollary 5.15]{IV:cofrob}. 
\end{proof}

A double application of \thref{charlQF} leads now immediately to the following characterization of quasi-Frobenius corings.

\begin{theorem}\thlabel{charQF}
Let $\cc$ be an $A$-coring. The following statements are equivalent.
\begin{enumerate}[(i)]
\item[(0)] The Morita context \equref{contextC*C} is strict;
\item $\cc$ is similar to $\*c$ as $A$-$\*c$ bimodule (i.e. $\cc$ is a quasi-Frobenius coring);
\item $\cc$ is similar to $\cc^*$ as $\cc^*$-$A$ bimodule;
\item $\cc$ is a direct summand of $(\*c)^n$ and of $(\cc^*)^m$, respectively as $A$-$\*c$ bimodule and $\cc^*$-$A$ bimodule, for certain numbers $n$ and $m$ (i.e.\ $\cc$ is left and right quasi-Frobenius);
\item $(\Ff,\Gg)$ is a left and right quasi-Frobenius pair;
\item $(\Gg,\Hh)$ is a left and right quasi-Frobenius pair;
\item the forgetful functor $\Ff$ and the functor $\Hh$ are similar (i.e.\ the induction functor $\Gg$ is quasi-Frobenius);
\item the functors ${_{\cc^*}\Hom_A}(-,\cc^*)$ and ${_{\cc^*}\Hom_A}(-,\cc):{_{\cc^*}\Mm_A}\to\Mm_{\*c^*}$ are similar;
\item the functors ${_{\cc^*}\Hom_A}(\cc^*,-)$ and ${_{\cc^*}\Hom_A}(\cc,-):{_{\cc^*}\Mm_A}\to{_{\*c^*}\Mm}$ are similar; 
\item left hand versions of (iii)-(viii), i.e.\ replace in the statements $\Mm^\cc$ by ${^\cc\Mm}$, $\Mm_A$ by ${_A\Mm}$, $\Ff$ by $\Ff'$, $\Gg$ by $\Gg'$, $\Hh$ by $\Hh'$, $\cc^*$ by $\*c$ and $\*c$ by $\cc^*$.
\end{enumerate}
If these statements are satisfied, than $\cc$ is finitely generated and projective as a left and right $A$-module.
\end{theorem}

Consider a ringmorphism $\iota:B\to A$. We say that $\cc$ is \emph{left $B$-locally quasi-Frobenius} if there exists and index set $I$ and a $B$-$\*c$ bimodule morphism $j:\cc\to (\*c)^{I}$, such that for all $c_1,\ldots,c_n\in\cc$, there exists an element $\bar{z}=(z_\ell)_{\ell\in I}\in (\cc^B)^{(I)}$ satisfying $c_i=\bar{z}\cdot j(c_i)=\sum_\ell z_\ell j_\ell(c_i)$ for all $i=1,\ldots, n$. If $I$ can be choosen to be a singleton, than $\cc$ is called \emph{left $B$-locally Frobenius}. 
Every $A$-locally quasi-Frobenius coring is left quasi-co-Frobenius (i.e. there exists an $A$-$\*c$ bimodule monomorphism $j:\ \cc\to(\*c)^{I}$ for some index set $I$) and the converse is true provided that $A$ is a PF-ring.

Locally (quasi-)Frobenius corings that are locally projective modules over their base algebra where characterized and be given a functorial description in \cite[Theorem 5.3, Theorem 5.16]{IV:cofrob}. To show their close relationship with the notion of a (left) quasi-Frobenius coring, we provide a partially new characterization in the next Theorem.

\begin{theorem}\thlabel{charlocQF}
Let $\cc$ be an $A$-coring and consider the Morita context associated to $\cc$ and $\*c$ in ${_B\Mm_{\*c}}$,
\begin{equation}\eqlabel{contextBC*C}
\NN_B(\cc,\*c)=({_B\End_{\*c}}(\cc),{_B\End_{\*c}}(\*c),{_B\Hom_{\*c}}(\*c,\cc),{_B\Hom_{\*c}}(\cc,\*c),\circ,\bul).
\end{equation}
Then the following statements hold,
\begin{enumerate}[(i)]
\item $\cc$ is left $B$-locally Frobenius if and only if there exists a left invertible element for~$\circ$;
\item $\cc$ is left $B$-locally quasi-Frobenius if and only if $\circ$ is locally surjective;
\end{enumerate}
Under any these conditions, $\cc$ is weakly locally projective as a right $A$-module.
\end{theorem}
\begin{proof}
\ul{(i)}. Let $j\in{_B\Hom_{\*c}}(\cc,\*c)$ be the left invertible element for $\circ$. Then for every finite set of elements $c_i\in\cc$ there exists an element $\barj\in{_B\Hom_{\*c}}(\*c,\cc)$, such that $\barj\circ j(c_i)=c_i$. Observe that $\phi:{_B\Hom_{\*c}}(\*c,\cc)\to \cc^B, \phi(f)=f(\varepsilon)$ is an isomorphism (see \cite[lemma 4.14]{IV:cofrob}). Putting $z=\phi(\barj)=\barj(\varepsilon)$, we obtain $z\cdot j(c_i)=c_i$, hence $\cc$ is $B$-locally Frobenius. The converse is proven in the same way.\\
\ul{(ii)}. Suppose that $\circ$ is locally surjective, then we know that for every finitely generated $\dd\subset \cc$, there exist a finite indexset $I_\dd$ and maps $j^\ell_\dd\in{_B\Hom_{\*c}}(\cc,\*c)$, $\barj^\ell_\dd\in{_B\Hom_{\*c}}(\*c,\cc)$ with $\ell\in I_\dd$ such that $\barj^\ell_\dd\circ j^\ell_\dd(d)=d$ for all $d\in\dd$. 
Put $I=\cup_{\dd} I_\dd$, where the union ranges over all finitely generated submodules $\dd$ of $\cc$. Then the universal property of the product induces a unique $B$-$\*c$ bilinear map $j:\cc\to \*c^I$, such that $j_i=\pi_i\circ j$, where $i\in I$ and $\pi_i: \*c^I\to \*c$ is the projection on the $i$-th component. Furthermore, for any finitely generated $\dd\subset \cc$, the maps $\barj^\ell_\dd\in{_B\Hom_{\*c}}(\*c,\cc)$ define an element $\bar{z}=\sum_{\ell\in I_\dd} \barj^\ell_\dd(\varepsilon)\in \cc^{(I_\dd)}\subset \cc^{(I)}$. The explicit form of $\phi^{-1}$ implies that
$$d=\sum_{\ell\in I_\dd} \barj^\ell_\dd \circ j^\ell_\dd (d) = \bar{z}\cdot j(d),$$
for all $d\in\dd$, i.e.\ $\cc$ is $B$-locally quasi-Frobenius.
Conversely, suppose that $\cc$ is $B$-locally quasi-Frobenius, and let $j:\cc\to \*c^I$ be the Frobenius map. Take any finitely generated submodule $\dd$ of $\cc$, then we know that there exists an element $\bar{z}=(z_\ell)_{\ell\in I}\in(\cc^{(I)})^B$ such that $\bar{z}\cdot j(d)=d$ for all $d\in \dd$. Denote by $I'$ the finite subset of $I$ containing all indices of non-zero entries of $\bar{z}$. Define for all $\ell\in I'$ maps $j^\ell\in{_B\Hom_{\*c}}(\cc,\*c)$ by $j^\ell=\pi_\ell\circ j$, where $\pi_\ell: \*c^I\to \*c$ are the canonical projections. Similarly, we define $B$-$\*c$ bilinear maps $\barj^\ell : \*c\to \cc,\ \barj^\ell(f)=z_\ell \cdot f$. Then we find for all $d\in\dd$, 
$$\sum_{\ell\in I'}\barj^\ell\circ j^\ell(d)= \sum_{\ell\in I'} z_\ell \cdot (\pi_\ell\circ j)(d)=\bar{z}\cdot j(d)=d,$$ 
hence  $\circ$ is locally surjective.

The last statement follows from \cite[Theorem 5.3(c)]{IV:cofrob}.
\end{proof}

From \thref{charlQF}, \thref{charQF} and \thref{charlocQF} it follows now immediately that a left quasi-Frobenius coring is left locally quasi-Frobenius and in particular a quasi-Frobenius coring is both left and right locally quasi-Frobenius.

\section{Locally projective Galois comodules}\selabel{Galois}

Let $\cc$ be an $A$-coring, $R$ a firm ring and $\Sigma\in{_R\Mm^\cc}$, then there is a ringmorphism $\iota\to\End^\cc(\Sigma)$ and we call $R$ the \emph{ring of coinvariants} of $\Sigma$. The Galois theory for $\Sigma$ studies the properties of the functor
$$-\ot_R\Sigma:\Mm_R\to\Mm^\cc$$
and its right adjoint
$$\Hom^\cc(\Sigma,-)\ot_RR:\Mm^\cc\to\Mm_R.$$
In this section we discuss the theory of Galois comodules whose ring of coinvariants is a ring with local units. This is a special situation of the theory developed in \cite{GTV:firm}, where the ring of coinvariants was supposed to be firm, so most result can be obtained from there. However some aspects of the theory differ slightly in this situation and stronger results can be obtained. In particular, we can characterize categories of comodules with a locally projective generator.

\begin{lemma}\lelabel{locprojgen}
Let $R$ be a ring with left local units, $\cc$ an $A$-coring that is flat as a left $A$-module and $\Sigma\in{_R\Mm^\cc}$. 
If the functor $-\ot_R\Sigma:\Mm_R\to \Mm^\cc$ establishes an equivalence of categories, then 
$\Sigma$ is a weakly locally projective in $\Mm^\cc$.
\end{lemma}

\begin{proof}
Denote $F=-\ot_R\Sigma:\Mm_R\to \Mm^\cc$ and let $G=\Hom^\cc(\Sigma,-)\ot_RR$ be the right adjoint of $F$.
Observe that $\Sigma\cong F(R)$, and therefore $R\cong GF(R)\cong G(\Sigma)$.
Consider a diagram of the form \equref{weaklyproj} in $\Mm^\cc$, with $X=\Sigma$. 
Since $E\cong FG(E)$, \leref{fg} implies that $G(E)$ is finitely generated in $\Mm_R$.
Apply the functor $G$ to the diagram \equref{weaklyproj}, then we obtain
\[
\xymatrix{
0 \ar[rr] && G(E) \ar[rr]^-{G(i)} && G(\Sigma)\cong R \ar[d]^{G(f)}\\
&& G(M) \ar[rr]^-{G(g)} && G(N) \ar[rr] && 0  
}
\]
Since $R$ is locally projective in $\Mm_R$ (see \leref{locprojcom}) we find a morphism $h:R\to G(M)$ making this diagram commutative on the image of $G(i)$. If we apply the functor $F$ to this diagram, then we see that $F(h):\Sigma\to M$ satisfies $g\circ F(h)\circ i=f\circ i$, i.e. $\Sigma$ is weakly locally projective in $\Mm^\cc$.
\end{proof}

\begin{theorem}\thlabel{locunits}
Let 
$\cc$ be an $A$-coring that is weakly locally projective as a left $A$-module and $\Sigma$ a right $\cc$-comodule. Then the following statements are equivalent
\begin{enumerate}[(i)]
\item There is a ring with left local units $R$ together with a ring morphism $\iota:R\to \End^\cc(\Sigma)$ 
and the functor $-\ot_R\Sigma:\Mm_R\to\Mm^\cc$ is an equivalence of categories;
\item there is a ring with left local units $R$ together with a ring morphism $\iota:R\to \End^\cc(\Sigma)$, 
$\Sigma$ is a weakly $R$-locally projective right $A$-module, 
$\can:\Sigma^*\ot_R\Sigma\to \cc$ is bijective and $\Sigma$ is faithfully flat as a left $R$-module; 
\item $\Sigma$ is a weakly locally projective generator in $\Mm^\cc$;
\end{enumerate}
Moreover, the ring $R$ is a two-sided ideal in $\End^\cc(\Sigma)$.
\end{theorem}

\begin{proof}
$\ul{(i)\Rightarrow(iii)}$.
We know from \cite[Theorem 5.9]{GTV:firm}, \cite[Theorem 3.4]{Ver:equiv} that $\Sigma$ is a generator in $\Mm^\cc$. It follows by \leref{locprojgen} that $\Sigma$ is weakly locally projective in $\Mm^\cc$.\\
${(iii)\Rightarrow(ii)}$.
Let us first construct the ring with local units $R$. Since $\cc$ is weakly locally projective as a left $A$-module, every finite subset of $\Sigma$ is contained in a subcomodule $F$ that is finitely generated as a right $A$-module (see \cite[19.12]{BrzWis:book}) and by \leref{locprojcom}, $F$ is finitely generated in $\Mm^\cc$. Hence $\Sigma$ is generated by its finitely generated subcomodules in $\Mm^\cc$. Denote  $\Omega=\bigoplus_{F\subset \Sigma} F$, the direct sum in $\Mm^\cc$ of all finitely generated subcomodules of $\Sigma$. Then there is a canonical surjection $\pi:\Omega\to \Sigma$ in $\Mm^\cc$. Let $R=P\ot_T Q$, where $P\subset \Hom^\cc(\Omega,\Sigma)$, consisting of all maps with finite support, that is, all $f\in \Hom^\cc(\Omega,\Sigma)$, such that $f$ is zero everywhere exept on a finitely generated direct summand of $\Omega$ (hence a finitely generated submodule of $\Sigma$), 
$Q=\Hom^\cc(\Sigma,\Omega)$ and $T={\End^\cc(\Omega)}$. For $f,g\in P$ and $\phi,\psi\in Q$ we define
$$(f\ot_T \phi)(g\ot_T \psi)=f\circ\phi\circ g\ot_T \psi = f\ot_T \phi\circ g\circ \psi.$$
By construction, $R$ is a two-sided ideal in $\End^\cc(\Sigma)$.
We claim that $R$ is a ring with left local units and $\Sigma$ is weakly $R$-locally projective. 
Take any $r=f\ot\phi\in R$, where $f\in P$, then $F=\im (f)$ is a subcomodule of $\Sigma$ that is finitely generated as a right $A$-module. By the local projectivity of $\Sigma$ in $\Mm^\cc$, we find a $\cc$-colinear map $\psi:\Sigma\to \Omega$ such that the following diagram commutes on the image of the inclusion $\iota$,
\[
\xymatrix{
0\ar[rr] && F\ar[rr]^-\iota && \Sigma \ar[lld]_-{\psi} \ar@{=}[d]\\
&& \Omega\ar[rr]_-\pi && \Sigma \ar[rr] && 0
}
\]
Since $F$ is finitely generated, $\psi\circ\iota(F)\subset \Omega$ is also finitely generated, and therefore contained in a direct summand $E$ of $\Omega$. Write $\Omega=E\oplus E'$ and define $\pi_E:\Omega\to \Sigma$, as $\pi|_E$ on $E$ and zero on $E'$. Then we find
$\pi_E\circ\psi\circ\iota =\pi\circ \psi\circ\iota =\iota$, hence $\pi_E\circ \psi\circ f= f$ and therefore 
$(\pi_E\ot \psi)\in R$ is a left local unit for $(f\circ \phi)$. 

A similar arguement shows that $R$ acts with local units on $\Sigma$. Furthermore, a morphism $R\to\Sigma\ot_A\Sigma^*$ is constructed as follows: since $\Sigma$ is weakly locally projective in $\Mm^\cc$, we know by \leref{locprojcom} that $\Sigma$ is weakly locally projective in $\Mm_A$. Hence for any element $f\ot \phi\in R$, we can find a local dual basis $\sum_i e_i\ot_A f_i\in \Sigma\ot_A \Sigma^*$ for the finitely generated $A$-module $\im f$. Then define $\iota(r)=\sum_i f\circ \phi(e_i)\ot_A f_i$. It follows that $\Sigma$ is weakly $R$-locally projective as a right $A$-module.

By the weak Galois structure theorem (see \cite[Theorem 5.9]{GTV:firm}), the generator property implies that $\can$ is an isomorphism and $\Sigma$ is flat as a left $R$-module. So we are finished if we prove the total faithfulness of $\Sigma$ as a left $R$-module.
Since $R$ is flat as left $R$-module (see \leref{locprojcom}), by \prref{faithfullyflatfirm}, it is enough to prove that $J\Sigma\neq \Sigma$ for any proper ideal $J\subset R$.
Arguments similar to the ones in \cite[18.4 (3)]{Wis:book} show that for any right ideal $J$ of $R$, the injective map
$J\ot_RR\to \Hom^\cc(\Sigma,J\Sigma)\ot_RR,\ j\ot_Rr\mapsto (x\mapsto j(x))\ot_Rr$ is an isomorphism. Details are as follows. Take $g\ot_R r\in  \Hom^\cc(\Sigma,J\Sigma)\ot_RR$. Then $g\ot_R r=ge\ot_R r$, where $e=\psi_e\ot_R \phi_e\in R$ is a left local unit for $r$. By definition of $R$, $\psi_e$ is only non-zero on a finitely generated submodule $E$ of $\Sigma$.
Let $\{u_1,\cdots,u_k\}$ be a set of generators of $E$ as right $A$-module, and complete it to a set of generators $\{u_i\}_{i\in I}$ for $\Sigma$.
Since $g(u_i)\in J\Sigma$, we can write $g(u_i)=\sum_jf_{ij}(u_j)$, with $f_{ij}\in J$. Let $J'$ be the subideal of
$J$ generated by the $f_{ij}$. Since  the $\{f_{ij}(u_i)\}$
generate $\im(g)$ as a right $A$-module, we have that $\im(g)\subset
J'\Sigma$. We relabel the generators $f_{ij}$ of $J'$ as 
$\{f_\ell\}_{\ell\in L}$.
Let $\pi_\ell:\ \Sigma^{(L)}\to \Sigma$ and $\iota_\ell:\ \Sigma\to \Sigma^{(L)}$, $\ell\in L$, be the natural projections and inclusions. Consider the following diagram
\[
\xymatrix{
0\ar[rr] && E \ar[rr]^-\iota && \Sigma \ar[d]^g\\
&& \Sigma^{(L)} \ar[rr]_-f && J'\Sigma \ar[rr] && 0
}
\]
where the map $f=\sum_{\ell\in L} f_\ell\circ\pi_\ell:\ \Sigma^{(L)}\to J'\Sigma$
is well defined and surjective. Since $\Sigma\in \Mm^\cc$ is weakly locally projective, there exists a map
$h:\Sigma\to \Sigma^{(L)}$ such that $g\circ \iota =f\circ h\circ\iota $, hence $g e= f\circ h e$. Using the universal property of the direct sum $\Sigma^{(L)}$, we can write $h=\sum_{\ell\in L} \iota_\ell\circ h_\ell:\ \Sigma\to \Sigma^{(L)}$ with $h_j\in \End^\cc(\Sigma)$.
Hence we obtain
$$f\circ h=\sum_{\ell, k\in L} f_k\circ\pi_k\circ \iota_\ell\circ h_\ell=\sum_{\ell\in L} f_\ell\circ h_\ell \in J'\subset J$$
where we used that since $J'$ is a right ideal of $R$ and $R$ is a right ideal of $\End^\cc(\Sigma)$, $J'$ is a right ideal of $\End^\cc(\Sigma)$. So we obtain that 
$$g\ot_R r= ge\ot_R r= f\circ he\ot_R r=f\circ h\ot_R r\in J\ot_R R,$$
proving $J\ot_RR \cong \Hom^\cc(\Sigma,J\Sigma)\ot_RR$.
If $J\neq R$, then we also have $\Hom^\cc(\Sigma,J\Sigma)\ot_RR\neq \Hom^\cc(\Sigma,R\Sigma)\ot_RR\cong R$ (since $R$ is an ideal in $\End^\cc(\Sigma)$, see \cite[Lemma 5.10]{GTV:firm}),
hence $J\Sigma\neq \Sigma$.\\
$\ul{(ii)\Rightarrow(i)}$.
Since $\Sigma$ is weakly $R$-locally projective, the ring $R$ acts with local units on $\Sigma$ with action induced by the ringmorphism $\iota:R\to\Sigma\ot_A\Sigma^*$, hence $\Sigma$ is firm as a left $R$-module and $\Sigma$ is also $R$-firmly projective in the sense of \cite{Ver:equiv}. Therefore, the implication follows from the Galois comodule structure theorem over firm rings, \cite[Theorem 3.4]{Ver:equiv}, \cite[Theorem 5.15]{GTV:firm}.\\
\end{proof}

Recall from \cite{BohmVer:cleft}, that to any $A$-coring $\cc$ and a right $\cc$-comodule $\Sigma$, we can associate a Morita context
\begin{equation}\eqlabel{contextSigma}
\CC(\Sigma)=(T,\*c,\Sigma,Q,\newtriangle,\newblacktriangle).
\end{equation}
Here we defined $T=\End^\cc(\Sigma)$, 
\begin{equation}
Q =\{\ q\in \Hom_A(\Sigma,{}^*\cc)\ |\ \forall x\in \Sigma,c\in \cc\quad
q(x_{[0]})(c)x_{[1]}=c_{(1)}q(x)(c_{(2)})\ \}
\end{equation}
with bimodule structure
\begin{eqnarray*}
(fq)(x)&= f q(x),\qquad &\text{for}\ f\in {}^*\cc,\ q\in Q,\ x\in
  \Sigma\qquad \text{and}\\
(qt)(x)&=q\big(t(x)\big),\qquad &\text{for}\ q\in Q,\ \ t\in T,\,\
  x\in \Sigma. 
\end{eqnarray*}
The bimodule structure on $\Sigma$ is given by
$$txf=t\big(x_{[0]}f(x_{[1]})\big),$$ 
for $t\in T$, $x\in \Sigma$ and $f\in {}^*\cc$.
And the Morita maps are
\begin{eqnarray}
\newblacktriangle:Q\ot_T \Sigma \to {}^*\cc,&\qquad& 
q\newblacktriangled x=q(x),\label{eq:F}\eqlabel{defF}\\
\newtriangle:\Sigma\ot_{{}^*\cc} Q\to T,&\qquad& 
x\newtriangled q(-) = x_{[0]}q(-)(x_{[1]}).\label{eq:G}\eqlabel{defG}
\end{eqnarray}
Recall from \cite[Remark 2.2]{BohmVer:cleft} that if $\cc$ is weakly locally projective as a left $A$-module, then $Q\cong \Hom_{\*c}(\Sigma,\*c)=\Hom^\cc(\Sigma,\Rat(\*c))$.

We say that a Morita context $\MM=(A,B,P,Q,\mu,\tau)$ can be \emph{restricted} to a (non-unital) subring $R\subset A$ if $\im\mu \subset R$. This restricted context is then constructed by considering $P$ and $Q$ in a canonical way as a left and right $R$-module and defining the map $\bar{\mu}:P\ot_B Q\to R$ as the coristriction of $\mu$, and $\bar{\tau}:P\ot_R Q\to B$ by
$$\xymatrix{\bar{\tau}: P\ot_R  Q\ar[rr]^-\pi  && P\ot_A Q \ar[rr]^-\tau && B},$$
where $\pi$ is the canonical projection.

If $\cc$ is weakly locally projective as a left $A$-module, then $Q=\Hom_{\*c}(\Sigma,\*c)\cong \Hom^\cc(\Sigma,\Rat(\*c))$, hence $q\newblacktriangled x=q(x)\in \Rat(\*c)$, for all $q\in Q$, so we can restrict the Morita context $\CC(\Sigma)$ to $\Rat(\*c)\subset \*c$.

If $R$ is a firm ring, $\iota:R\to \End^\cc(\Sigma)$ a ring morphism such that $R$ is a right ideal of $T$ and $\Sigma$ becomes a firm left $R$-module (e.g.\ $\Sigma$ is a weakly locally projective generator in $\Mm^\cc$ and $R$ is as constructed in \thref{locunits}) then 
for all $q\in Q$ and $x\in \Sigma$, 
$$x \newtriangled q= r\cdot x^r_{[0]}q(-)(x^r_{[1]})\in R$$
where we used that $R$ is a right ideal of $T$.
Therefore, the Morita context $\CC(\Sigma)$ can be restricted to $R\subset T$.

If it exists, we denote the (twofold) restricted Morita context as follows, 
\begin{equation}\eqlabel{restricted}
\ol{\cc}(\Sigma)=(R,\Rat(\*c),\Sigma,Q,\bar{\newtriangle},\bar{\newblacktriangle}).
\end{equation}

If $\Rat(\*c)$ is dense in the finite topology on $\*c$ (this last condition is for example satisfied if $\cc$ is a right $B$-locally quasi-Frobenius coring \cite[Propositon 5.9]{IV:cofrob}), then we know (see \cite[Lemma 4.13]{CDV:comatrix}) that there is an isomorphism of categories
$$\Mm^\cc \to \Mm_{\Rat(\*c)}.$$

The previous observations lead to a further equivalent statement for the conditions of \thref{locunits}. First we prove the following lemma about Morita contexts over non-unital rings. 

\begin{lemma}\lelabel{contextlocunit}
Let $(R,S,P,Q,\mu,\tau)$ be a Morita context, where $R$ is a ring with local units, and $R$ acts with local units on $P$. Then the surjectivity of $\mu$ implies its bijectivity.
\end{lemma}
\begin{proof}
Take any $\sum p\ot_Sq\in\ker\mu$ and let $e\in R$ be a local unit for $p$. By the surjectivity of $\mu$ we can write $e=\mu(p_e\ot_Sq_e)$, hence
$$p\ot_Sq=ep\ot_Sq=\mu(p_e\ot_Sq_e)p\ot_Sq=p_e\tau(q_e\ot_Rp)\ot_Sq=p_e\ot_Sq_e\mu(p\ot_Sq)=0.$$
\end{proof}

\begin{theorem}\thlabel{strstr2}
Let $\cc$ be an $A$-coring that is weakly locally projective as a left $A$-module, and suppose that $\Rat(\*c)$ is dense in the finite topology on $\*c$. For a right $\cc$-comodule $\Sigma$, the following statements are equivalent
\begin{enumerate}[(i)]
\item $\Sigma$ is a weakly locally projective generator in $\Mm^\cc$;
\item there is a ring with local units $R$, together with a ringmorphism $\iota:R\to \End^\cc(\Sigma)$, such that $R$ is a right ideal in $\End^\cc(\Sigma)$, $R$ acts with local units on $\Sigma$ and the maps of the Morita context \equref{restricted} are surjective.
\end{enumerate}
\end{theorem}
\begin{proof}
It follows form \leref{contextlocunit} that the surjectivity of the Morita maps implies their bijectivity. Hence the strict Morita context $\ol{\CC}(\Sigma)$ induces an equivalence between the categories $\Mm_R$ and $\Mm_{\Rat(\*c)}\cong \Mm^\cc$. This equivalence is given by the functor $-\ot_R\Sigma$. By \thref{locunits} this statement is equivalent with statement (i).
\end{proof}

\begin{theorem}\thlabel{weakstrSigma}
Let $\cc$ be an $A$-coring that is weakly locally projective as a left $A$-module, and suppose that $\Rat(\*c)$ is dense in the finite topology on $\*c$. Let $R$ be a firm ring and suppose that $\Sigma\in{_R\Mm^\cc}$, such that $R$ is a right ideal in $\End^\cc(\Sigma)$.
Then the following statements are equivalent
\begin{enumerate}[(i)]
\item $\Sigma$ is a generator in $\Mm^\cc$;
\item $\Sigma$ is flat as a left $R$-module and a $R$-comonadic Galois comodule;
\item the functor $\Hom^\cc(\Sigma,-)\ot_RR:\Mm^\cc\to\Mm_R$ is fully faithful;
\item the map $\bar{\newblacktriangle}$ of the Morita context \equref{restricted} is surjective;
\end{enumerate}
\end{theorem}
\begin{proof}
The equivalence of (i)-(ii)-(iii) follows from the weak sturcture theorem for Galois comodules (see \cite[Theorem 5.9]{GTV:firm}). 
So, we only have to prove the equivalence between (iii) and (iv). 
The functor $\Hom^\cc(\Sigma,-)\ot_RR$ will be fully faithfull if and only if the counit 
$$\zeta_M:\Hom^\cc(\Sigma,M)\ot_R\Sigma\to M,\ \zeta_M(\phi\ot_R x)=\phi(x)$$
of the adjunction $(-\ot_R\Sigma,\Hom^\cc(\Sigma,-)\ot_RR)$ is a natural isomorphism. 
Remark that since $Q\cong \Hom^\cc(\Sigma,\Rat(\*c))$, we have $\zeta_{\Rat(\*c)}=\bar{\newblacktriangle}$, so if $\Hom^\cc(\Sigma,-)\ot_RR$ is fully faithful, then $\bar{\newblacktriangle}$ is bijective. Conversely, an inverse for $\zeta_M$ is constructed as follows. Since $\Rat(\*c)$ is dense in the finite topology of $\*c$, for all $M\in\Mm^\cc$ and all $m\in M$, there exists an element $e\in\Rat(\*c)$ such that $m\cdot e=m$. Now take $q_e\ot_R x_e\in Q\ot_R \Sigma$ such that $q_e\bar{\newblacktriangled} x_e= e$. We define 
$$\theta_M:M\to \Hom^\cc(\Sigma,M)\ot_R\Sigma,\ \theta_M(m)= mq_e(-)\ot_R x_e.$$
One can check that $\theta_M$ is well defined and natural in $M$. Moreover, $\zeta_M\circ\theta_M (m)=m\cdot q_e(x_e)=m\cdot e=m$ and $\theta_M\circ \zeta_M(\phi\ot_R x)=\phi(x)q_e(-)\ot_R x_e$,
where $e$ can be choosen to be a right local unit for $x$. Furthermore, using the fact that $xq_e(-)\in T$ in combination with \cite[Lemma 5.11]{GTV:firm}, we find $\phi(x)q_e(-)\ot_R x_e=\phi\ot_R xq_e(x_e)=\phi\ot_R x\cdot e=\phi\ot_R x$, hence $\theta_M$ is a two-sided inverse for $\zeta_M$.
\end{proof}

\section{The coring as Galois comodule}\lelabel{CisGalois}

In this section we show how the theory of locally quasi-Frobenius corings is related to the theory of Galois comodules with firm coinvariant rings. We show that locally quasi-Frobenius corings are precizely corings that are a weakly locally projective generator in the category of their left or right comodules, i.e.\ faithfully flat \emph{infinite} Galois comodules. As a consequence, quasi-Frobenius corings are progenerators in the category of their left or right comodules, i.e.\ faithfully flat \emph{finite} Galois comodules.

Let $\cc$ be an $A$-coring that is weakly locally projective as a right $A$-module. Denote $R=\Rat(\cc^*)$ and consider the pair of adjoint functors
\begin{equation}\eqlabel{ccGalois}
\xymatrix{ {\Mm_R}
\ar@<0.5ex>[rrr]^-{-\ot_R\cc} &&&
{\Mm^\cc\ .}
\ar@<0.5ex>[lll]^{{\Hom^\cc}(\cc,-)\ot_RR }
}
\end{equation}
We want to examine when this pair of adjoint functors is an equivalence of categories. The comatrix coring is given by ${\Hom^\cc}(\cc,\cc)\ot_RR\ot_R\cc \cong\cc^*\ot_R\cc$, hence we can calculate the canonical map as follows
$$\can:\cc^*\ot_R\cc\cong \cc,\ \can(f\ot_Rc)=f(c_{(1)})c_{(2)}=f\cdot c.$$
Because of the counit property, $\can$ is clearly surjective. If moreover $R$ is dense in the finite topology on $\cc^*$, then, by the fact that $R$ is an ideal in $\cc^*$ and \cite[Lemma 5.11]{GTV:firm}, we obtain that $\can$ is bijective with inverse 
$$\can^{-1}:\cc\to \cc^*\ot_R \cc,\ \can^{-1}(c)=\varepsilon\ot_R c.$$ 

Applying the weak and strong structure theorem for Galois comodules (see \cite[Theorem 5.9 and Theorem 5.15]{GTV:firm}), we obtain immediately the following.
\begin{theorem}\thlabel{strstrC}
If $\cc$ is an $A$-coring that is flat as a left $A$-module and weakly locally projective as a right $A$-module such that $R=\Rat(\cc^*)$ is dense in the finite topology on $\cc^*$,
then 
\begin{enumerate}
\item the functor $\Hom^\cc(\cc,-)\ot_RR$ is fully faithful if and only if any of the following equivalent conditions holds
\begin{enumerate}[(i)]
\item $\cc$ is a generator in $\Mm^\cc$;
\item $\cc$ is flat as a left $R$-module;
\item the functor $\Hom^\cc(\cc,-):\Mm^\cc\to \Mm_{\cc^*}$ is fully faithful.
\end{enumerate}
\item $(-\ot_R\cc,{\Hom^\cc}(\cc,-)\ot_RR)$ is a pair of inverse equivalences between the categories ${\Mm_R}$ and ${\Mm^\cc}$ if and only if any of the following equivalent conditions holds
\begin{enumerate}[(i)]
\item $\cc$ is faithfully flat as a left $R$-module;
\item $\cc$ is a generator in $\Mm^\cc$ such that $-\ot_R\cc:\Mm_R\to\Mm^\cc$ is faithful;
\item $\cc$ is flat as a left $R$-module and $\Rat(\cc^*)$ is coflat as a left $\cc$-comodule.
\end{enumerate}
\end{enumerate}
\end{theorem}

\begin{remark}
Statement (2)(iii) of \thref{strstrC} follows in fact from a generalization of the strong structure theorem for Galois comodules, given in \cite[Theorem 4.27]{V:PhD}.
\end{remark}

To study the Galois comodule $\cc$ we can consider the Morita context \equref{contextSigma}, which is in this case,
\begin{equation}\eqlabel{contextC}
\CC(\cc)=(\End^\cc(\cc)\cong\cc^*,\*c,\cc,Q,\newtriangle,\newblacktriangle),
\end{equation}
If we suppose again that $\cc$ is weakly locally projective as a left $A$-module, then $Q=\Hom_{\*c}(\cc,\Rat(\*c))\cong \Hom_{\*c}(\cc,\*c)\cong {_{\cc^*}}\Hom(\cc,\cc^*)$, where the last isomorphism is given by switching the arguements (see \cite[Proposition 4.2]{IV:cofrob}). If moreover $\cc$ is weakly locally projective as a  right $A$-module, then 
$Q\cong {_{\cc^*}\Hom}(\cc,\Rat(\cc^*))=:\tilde{Q},$
and we can consider the restriced Morita context
\begin{equation}\eqlabel{contextRS}
\ol{\CC}(\cc)=(S:=\Rat(\cc^*),R:=\Rat(\*c),\cc,Q,\bar{\newtriangle},\bar{\newblacktriangle}),
\end{equation}
where 
\begin{eqnarray*}
\bar{\newtriangle}:\ \cc\otimes_{\*c}\tilde{Q}\to\Rat(\cc^*),&&c\bar{\newtriangled} \tilde{q}=\tilde{q}(c);\\
\bar{\newblacktriangle}:\ Q\otimes_{\cc^*}\cc\to\Rat(\*c),&& q\bar{\newblacktriangled} c=q(c).
\end{eqnarray*}

\begin{lemma}\lelabel{isocontexts}
Let $\cc$ be an $A$-coring that is weakly locally projective as a left $A$-module. Then there are isomorphisms of Morita contexts 
$$\NN_k(\cc,\*c) \cong \NN_k^{\rm top}(\cc,\cc^*)\cong \CC(\cc)$$ 
where the first two contexts are constructed as in \equref{contextBC*C} and the last context is the one of \equref{contextC}.
\end{lemma}

\begin{proof}
The isomorphism between the first two contexts is proven as in \cite[Theorem 5.10]{IV:cofrob}.

For the last isomorphism, just observe that $\End_{\*c}(\cc)\cong \End^\cc(\cc)\cong \cc^*$,
$\End_{\*c}(\*c)\cong\*c$,
$\Hom_{\*c}(\*c,\cc)\cong \cc$ and 
$\Hom_{\*c}(\cc,\*c)\cong Q$. One can easily check that these isomorphisms induce an isomorphism of Morita contexts between $\NN_k(\cc,\*c)$ and $\CC(\cc)$.

Similarly,
${_{\cc^*}\End}(\cc)\cong {^\cc\End}(\cc)\cong {\*c}$,
${_{\cc^*}\End}(\cc^*)\cong\cc^*$,
${_{\cc^*}\Hom}(\cc^*,\cc)\cong \cc$ and 
${_{\cc^*}\Hom}(\cc,\cc^*)\cong \tilde{Q}$, what induces the isomorphism $\NN_k^{\rm top}(\cc,\cc^*)\cong \CC(\cc)$.
\end{proof}

\begin{theorem}\thlabel{structureQF1}
Let $\cc$ be an $A$-coring that is weakly locally projective as a left $A$-module. 
The following statements are equivalent.
\begin{enumerate}[(i)]
\item $\cc$ is left locally $k$-locally quasi-Frobenius;
\item $\cc$ is weakly locally projective as a right $A$-module, $R=\Rat(\cc^*)$ is dense in the finite topology on $\cc^*$ and 
the map $\bar{\newtriangle}$ of the Morita context \equref{contextRS} is surjective (hence bijective).
\item the functor $-\ot_R\cc: \Mm_R\to \Mm_{\Rat(\*c)}$ is fully faithful;
\end{enumerate}
If moreover $S=\Rat(\*c)$ is dense in the finite topology on $\*c$, the previous statements are furthermore equivalent to any of the following assertions,
\begin{enumerate}[(i)]
\item[(iv)] the functor $-\ot_R\cc: \Mm_R\to \Mm^\cc$ is fully faithful;
\item[(v)] $\cc$ is a generator in ${^\cc\Mm}$;
\item[(vi)] $\cc$ is flat as a right $S$-module;
\item[(vii)] the functor ${^\cc\Hom}(\cc,-):{^\cc\Mm}\to {_{\*c}\Mm}$ is fully faithful;
\item[(viii)] the functor $S\ot_S{^\cc\Hom}(\cc,-):{^\cc\Mm}\to {_S\Mm}$ is fully faithful.
\end{enumerate}
\end{theorem}

\begin{proof}
$\ul{(i)\Leftrightarrow (ii)}$.
By \thref{charlocQF}, we know that $\cc$ is left $k$-locally quasi-Frobenius if and only if the Morita map $\circ$ of 
$$\NN(\cc,\*c)=({\End_{\*c}}(\cc),{\End_{\*c}}(\*c),{\Hom_{\*c}}(\*c,\cc),{\Hom_{\*c}}(\cc,\*c),\circ,\bul)$$ 
is locally surjective. 
The isomorphisms of Morita contexts in \leref{isocontexts} imply the following algebra and bimodule maps
\begin{eqnarray}
\psi:{\End_{\*c}}(\cc)\cong \cc^* && \psi(f)=\varepsilon\circ f,\eqlabel{psi}\\
\phi:{\Hom_{\*c}}(\cc^*,\cc)\to \cc && \phi(j)=j(\varepsilon), \\
\bar{\phi}: {\Hom_{\*c}}(\cc,\*c)\to \tilde{Q} && \bar{\phi}(\barj)(c)(d)= \barj(d)(c). \eqlabel{barphi}
\end{eqnarray}
Since these constitute a morphism of Morita context we have moreover that
\begin{equation}\eqlabel{jbarjcomp}
\psi( j\circ \barj)= \phi(j)\newtriangle \bar{\phi}(\barj).
\end{equation}
Take any element $f\in\Rat(\cc^*)$, and a finite number of representants $c_i\in \cc$ and $f_i\in\Rat(\cc^*)$ such that $f_{[-1]}\ot_A f_{[0]}= \sum_i c_i\ot_A f_i\in \cc\ot_A\Rat(\cc^*)$. By the local surjectivity of $\circ$, there exist elements $j_\ell\in {\Hom_{\*c}}(\*c,\cc)$ and $\barj_\ell\in {\Hom_{\*c}}(\cc,\*c)$ such that $\sum_\ell j_\ell\circ \barj_\ell(c_i)= c_i$ for all $c_i$. 
Denote $\phi(j_\ell)=j_\ell(\varepsilon)= z_\ell$ and $\bar{\phi}(\barj_\ell) =\tilde{q}_\ell$. Then we have because of \equref{jbarjcomp} 
\begin{eqnarray*}
\varepsilon(c_i)&=& \sum_\ell \varepsilon( j_\ell\circ \barj_\ell(c_i)) 
=\psi( j_\ell\circ \barj_\ell)(c_i)\\
&=& \sum_\ell (z_\ell \newtriangle \tilde{q}_\ell)(c_i)=\sum_\ell \tilde{q}_\ell(z_\ell)(c_i).
\end{eqnarray*} 
Hence we find
\begin{eqnarray*}
f 
&=& \sum_i  \varepsilon(c_i)f_i
= \sum_{i, \ell} (z_\ell \newtriangle \tilde{q}_\ell)(c_i) f_i \\
&=& (\sum_\ell z_\ell \newtriangle \tilde{q}_\ell)\cdot f 
= (\sum_\ell z_\ell \newtriangle \tilde{q}_\ell f) 
\end{eqnarray*}
where we used in the last equality that $\newtriangle$ is a right $\cc^*$-linear map.
So $\newtriangle$ is surjective onto $\Rat(\cc^*)$, i.e.\ $\bar{\newtriangle}$ is surjective.

Conversely, suppose that $\bar{\newtriangle}$ is surjective, than we have to show that $\circ$ is locally surjective. Take any finite set $c_i\in \cc$. Since $\Rat(\cc^*)$ is dense in the finite topology of $\cc^*$, we know that there exists a local unit $e\in\Rat(\cc^*)$ such that $e\cdot c_i=e(c_{i(1)}) c_{i(2)}=c_i$ for all $c_i$. By the surjectivity of $\bar{\newtriangle}$, we can write $e=z_e\newtriangle \tilde{q}_e$ for certain elements $\tilde{q}_e\in \tilde{Q}$ and $z_e\in \cc$. Using the isomorphisms \equref{psi}-\equref{barphi}, we obtain elements $j_\ell=\phi^{-1}(z_e)$ and $\barj_\ell=\bar{\phi}^{-1}(\tilde{q}_e)$ such that $\sum_\ell j_\ell\circ \barj_\ell(c_i)=\psi^{-1}(z_e\newtriangle \tilde{q}_e)(c_i)= (q_e\newblacktriangle z_e)\cdot c_i=e\cdot c_i=c_i$.

$\ul{(ii)\Leftrightarrow (iii)}$. Follows from Morita theory.

$\ul{(iii)\Leftrightarrow (iv)}$. Follows from the fact that $\Mm^\cc\cong\Mm_S$ if $S$ is dense in $\*c$.

To prove the equivalence with the other statements, observe that by left-right duality the construction of the Morita context \equref{contextSigma}, can be repeated for a left $\cc$-comodule. Consider the Morita context $\widetilde{\CC}(\cc)$, associated to $\cc$ as left $\cc$-comodule in this way, 
$$\widetilde{\CC}(\cc)=({^\cc\End}(\cc)\cong\*c^{\rm op}, (\cc^*)^{\rm op},\cc, \tilde{Q}, \tilde{\newtriangle}, \tilde{\newblacktriangle}).$$
One can verify that there is an isomorphism of Morita contexts $\widetilde{\CC}(\cc)^{\rm top}\cong \CC(\cc)$. Therefore, the surjectivity of $\bar{\newtriangle}$ implies that the restriction of $\tilde{\newblacktriangle}$ onto $\Rat(\cc^*)$ will be surjective. The equivalence with statements $(v)-(viii)$ follows then from (the left hand version of) \thref{weakstrSigma}.
\end{proof}

\begin{theorem}\thlabel{structureQF2}
Let $\cc$ be an $A$-coring.
The following statements are equivalent.
\begin{enumerate}
\item $\cc$ is left and right $k$-locally quasi-Frobenius;
\item $\cc$ is weakly locally projective as a left $A$-module $S=\Rat(\*c)$ is dense in the finite topology on $\*c$ and $\cc$ is a weakly locally projective generator in $\Mm^\cc$;
\item $\cc$ is weakly locally projective as a right $A$-module, $R=\Rat(\cc^*)$ is dense in the finite topology on $\cc^*$ and $\cc$ is a weakly locally projective generator in ${^\cc\Mm}$; 
\item $\cc$ is weakly locally projective as a left and right $A$-module, $R=\Rat(\cc^*)$ and $S=\Rat(\*c)$ are dense in the finite topology on respectively $\cc^*$ and $\*c$, and any of the following conditions holds
\begin{enumerate}[(i)]
\item the Morita context \equref{contextRS} is strict;
\item $-\ot_R\cc:\Mm_R\to \Mm^\cc$ is an equivalence of categories;
\item $\cc$ is faithfully flat as a left $R$-module;
\item $\cc\ot_S-:{_S\Mm}\to {^\cc\Mm}$ is an equivalence of categories;
\item $\cc$ is faithfully flat as a right $S$-module.
\end{enumerate}
\end{enumerate}
\end{theorem}

\begin{proof}
The equivalence of $(1)$ and $(4)(i)$ follows from a double application of \thref{structureQF1}.\\
Let us prove the equivalence between $(2)$ and $(4)(i)$.  By \thref{strstr2}, we know already that $(4)(i)$ implies $(2)$ and that $(2)$ implies the existence of a ring with local units $B$ that is an ideal in $\cc^*$ and such that $\newtriangle$ is surjective onto $B$. However, since $\im\newtriangle\subset\Rat(\cc^*)$, we find that $\Rat(\cc^*)$ is a as well a ring with local units, hence dense in the finite topology on $\cc^*$. This shows that $(2)$ implies $(4)(i)$. \\
The equivalence with $(4)(ii) - (4)(iv)$ follows now from \thref{locunits}, since the canonical map is always an isomorphism.\\
The equivalence between $(3)$ and $(4)$ follows by symmetric arguements.
\end{proof}

\begin{remark}
If $C$ is a coalgebra over a commutative ring, then (using the notation of \thref{structureQF2}) $C^*={^*C}$ and $R=S$. Hence $R$ is dense in the finite topology of $C^*$ if and only if $S$ is. A similar argumentation as in the proof of \thref{structureQF2} $(2)\Rightarrow (4)(i)$ shows that if $C$ is weakly locally projective as a $k$-module and $C$ is a weakly locally projective generator in $\Mm^C$, then $R$($=S$) will be dense in $C^*$. 
Therefore, we obtain from \thref{structureQF2} the following characterization.
\begin{enumerate}[(i)]
\item[] \hspace{-1cm} If $C$ is weakly locally projective as $k$-module, then the following are equivalent,
\item  $C$ is a $k$-locally quasi-Frobenius coalgebra;
\item $C$ is a weakly locally projective generator in $\Mm^C$;
\item $C$ is a weakly locally projective generator in ${^C\Mm}$.
\end{enumerate}
\end{remark}

We know that a right (respectively left) $A$-module $M$ is weakly $B$-locally projective, were $B$ is a ring with unit, if and only if $M$ is finitely generated and projective as a right (respectively left) $A$-module. Hence \thref{structureQF1} and \thref{structureQF2} yield immediately the following structure theorems for what we will term \emph{(right) $k$-quasi-Frobenius corings}.

\begin{corollary}\colabel{structureQF}
Let $\cc$ be an $A$-coring that is finitely generated and projective as a left $A$-module.
The following statements are equivalent.
\begin{enumerate}
\item $\cc$ is a direct summand of a finite number of copies of $\cc^*$ (i.e.\ $\cc$ is right $k$- quasi-Frobenius);
\item $\cc$ is finitely generated and projective as a right $A$-module, 
and any of the following conditions hold
\begin{enumerate}[(i)]
\item in the context \equref{contextC}, the map ${\newtriangle}$ is surjective (hence bijective);
\item the functor $-\ot_{\cc^*}\cc:\Mm_{\cc^*}\to \Mm^\cc$ is fully faithful;
\item $\cc$ is a generator in ${^\cc\Mm}$;
\item $\cc$ is flat as a right $\*c$-module;
\item the functor ${^\cc\Hom}(\cc,-):{^\cc\Mm}\to {_{\*c}\Mm}$ is fully faithful.
\end{enumerate}
\end{enumerate}
\end{corollary}

\begin{corollary}\colabel{structureQF2}
Let $\cc$ be an $A$-coring.
The following statements are equivalent.
\begin{enumerate}
\item $\cc$ and $\cc^*$ are similar as left $\cc^*$-modules (i.e.\ is $k$-quasi-Frobenius);
\item $\cc$ and $\*c$ are similar as right $\*c$-modules;
\item $\cc$ finitely generated and projective as a left $A$-module and $\cc$ is a finitely generated and projective generator in $\Mm^\cc$;
\item $\cc$ is finitely generated and projective projective as a right $A$-module and $\cc$ is a finitely generated and projective generator in ${^\cc\Mm}$; 
\item $\cc$ is finitely generated and projective as a left and right $A$-module and any of the following conditions hold
\begin{enumerate}[(i)]
\item the Morita context \equref{contextC} is strict, hence $-\ot_{\cc^*}\cc:\Mm_{\cc^*}\to \Mm_{\*c}$ is an equivalence of categories;
\item $-\ot_{\cc^*}\cc:\Mm_{\cc^*}\to \Mm^\cc$ is an equivalence of categories;
\item $\cc$ is faithfully flat as a left $\cc^*$-module;
\item $\cc\ot_{\*c}-:{_{\*c}\Mm}\to {^\cc\Mm}$ is an equivalence of categories;
\item $\cc$ is faithfully flat as a right ${\*c}$-module.
\end{enumerate}
\end{enumerate}
\end{corollary}


\begin{thebibliography}{10}
\expandafter\ifx\csname url\endcsname\relax
  \def\url#1{{\tt #1}}\fi
\expandafter\ifx\csname urlprefix\endcsname\relax\def\urlprefix{URL }\fi
\providecommand{\eprint}[2][]{\url{#2}}

\bibitem{BohmVer:cleft}
G.~B{\"o}hm and J.~Vercruysse, Morita theory for coring extensions and cleft
  bicomodules, {\em Adv. Math.\/}, {\bf 209}~(2), (2007) 611--648.

\bibitem{BrzWis:book}
T.~Brzezinski and R.~Wisbauer, Corings and comodules, volume 309 of {\em London
  Mathematical Society Lecture Note Series\/}, Cambridge University Press,
  Cambridge, 2003.

\bibitem{CDV:comatrix}
S.~Caenepeel, E.~De~Groot and J.~Vercruysse, Galois theory for comatrix
  corings: descent theory, {M}orita theory, {F}robenius and separability
  properties, {\em Trans. Amer. Math. Soc.\/}, {\bf 359}~(1), (2007) 185--226
  (electronic).

\bibitem{CMZ}
S.~Caenepeel, G.~Militaru and S.~Zhu, Frobenius and separable functors for
  generalized module categories and nonlinear equations, volume 1787 of {\em
  Lecture Notes in Mathematics\/}, Springer-Verlag, Berlin, 2002.

\bibitem{CIGTN:Frob}
F.~Casta{\~n}o~Iglesias, J.~G{\'o}mez-Torrecillas and
  C.~N{\u{a}}st{\u{a}}sescu, Frobenius functors: applications, {\em Comm.
  Algebra\/}, {\bf 27}~(10), (1999) 4879--4900.

\bibitem{IglNas:QF}
F.~Casta{\~n}o~Iglesias, C.~Nastasescu and J.~Vercruysse, Quasi-{F}robenius functors. Applications, preprint 2008.

\bibitem{EGT:comatrix}
L.~El~Kaoutit and J.~G{\'o}mez-Torrecillas, Comatrix corings: {G}alois corings,
  descent theory, and a structure theorem for cosemisimple corings, {\em Math.
  Z.\/}, {\bf 244}~(4), (2003) 887--906.

\bibitem{GTN}
J.~G{\'o}mez~Torrecillas and C.~N{\u{a}}st{\u{a}}sescu, Quasi-co-{F}robenius
  coalgebras, {\em J. Algebra\/}, {\bf 174}~(3), (1995) 909--923.

\bibitem{GTV:firm}
J.~G{\'o}mez-Torrecillas and J.~Vercruysse, Comatrix corings and {G}alois
  comodules over firm rings, {\em Algebr. Represent. Theory\/}, {\bf 10}~(3),
  (2007) 271--306.

\bibitem{Guo}
G.~Guo, Quasi-{F}robenius corings and quasi-{F}robenius extensions, {\em Comm.
  Algebra\/}, {\bf 34}~(6), (2006) 2269--2280.

\bibitem{IV:cofrob}
M.~Iovanov and J.~Vercruysse, Co-{F}robenius corings and adjoint functors, {\em
  J. Pure Appl. Algebra\/}, {\bf 212}~(9), (2008) 2027--2058.

\bibitem{KelStr:2cat}
G.~M. Kelly and R.~Street, Review of the elements of {$2$}-categories, in:
  Category Seminar (Proc. Sem., Sydney, 1972/1973), pp. 75--103. Lecture Notes
  in Math., Vol. 420, Springer, Berlin, 1974.

\bibitem{Lin:Sem}
B.~I.-p. Lin, Semiperfect coalgebras, {\em J. Algebra\/}, {\bf 49}~(2), (1977)
  357--373.

\bibitem{Mul:Morita}
B.~J. M{\"u}ller, The quotient category of a {M}orita context, {\em J.
  Algebra\/}, {\bf 28}, (1974) 389--407.

\bibitem{QFcodes}
A.~A. Nechaev, A.~S. Kuz'{}min and V.~T. Markov, Linear codes over finite rings
  and modules, {\em Fundam. Prikl. Mat.\/}, {\bf 3}~(1), (1997) 195--254,
  functional analysis, differential equations and their applications (Russian)
  (Puebla, 1995).

\bibitem{NicYou:QF}
W.~K. Nicholson and M.~F. Yousif, Quasi-{F}robenius rings, volume 158 of {\em
  Cambridge Tracts in Mathematics\/}, Cambridge University Press, Cambridge,
  2003.

\bibitem{Sw}
M.~Sweedler, The predual theorem to the {J}acobson-{B}ourbaki theorem, {\em
  Trans. Amer. Math. Soc.\/}, {\bf 213}, (1975) 391--406.

\bibitem{V}
J.~Vercruysse, Local units versus local projectivity dualisations: corings with
  local structure maps, {\em Comm. Algebra\/}, {\bf 34}~(6), (2006) 2079--2103.

\bibitem{V:PhD}
J.~Vercruysse, Galois theory for corings and comodules, 2007, {P}hD-thesis,
  Vrije Universiteit Brussel, \urlprefix\url{homepages.vub.ac.be/~jvercruy}.

\bibitem{Ver:equiv}
J.~Vercruysse, Equivalences between categories of modules and categories of
  comodules, {\em Acta Math. Sin.\/}, in press,
  \urlprefix\url{http://arxiv.org/math.RA/0604423}.

\bibitem{Wis:book}
R.~Wisbauer, Foundations of module and ring theory, volume~3 of {\em Algebra,
  Logic and Applications\/}, Gordon and Breach Science Publishers,
  Philadelphia, PA, 1991, german edition, a handbook for study and research.

\bibitem{ZH}
B.~Zimmermann-Huisgen, Pure submodules of direct products of free modules, {\em
  Math. Ann.\/}, {\bf 224}~(3), (1976) 233--245.

\end{thebibliography}
\end{document}